\documentclass[11pt]{article}

\usepackage{amsmath,amssymb,amsthm,mathtools}
\usepackage[left=2.0cm,right=2.0cm,top=2.3cm,bottom=2.3cm]{geometry}
\usepackage{enumitem}
\usepackage{cite}
\usepackage[colorlinks=true,citecolor=red,linkcolor=blue,urlcolor=blue]{hyperref}

\newtheorem{theorem}{Theorem}[section]
\newtheorem{corollary}[theorem]{Corollary}
\newtheorem{proposition}[theorem]{Proposition}
\newtheorem{lemma}[theorem]{Lemma}
\newtheorem{definition}[theorem]{Definition}
\newtheorem{remark}[theorem]{Remark}
\newtheorem{example}[theorem]{Example}
\numberwithin{equation}{section}

\DeclareMathOperator{\Ric}{Ric}

\newcommand{\RR}{\mathbb{R}}

\newcommand{\dd}{\,\mathrm{d}}
\newcommand{\Cone}{\mathcal{C}}

\title{Angular-Pole Flag Curvature and Conformal Control of Radial Jacobi Fields in Diagonal Conic Finsler Optics}

\author{
 Esra Sengelen Sevim\\
 Istanbul Ni\d{s}anta\d{s}ı University  \\
 \texttt{esra.sevim@nisantasi.edu.tr}
 \and
 Semail \"{U}lgen\thanks{Corresponding author: \texttt{sulgen@antalya.edu.tr}}\\
 Antalya Bilim University, Antalya, Turkey
}
\begin{document}
\maketitle

\begin{abstract}
We study diagonal conic Finsler optical structures with radial coefficients and an arbitrary angular Finsler norm. We obtain the spray curvature block and identify the angular-pole radial-transverse flag curvature at a critical radius with the negative optical-arclength radial Jacobi coefficient. We prove that a compactly supported radial conformal two-jet determines the critical radius and any signed radial coefficient. A navigation Randers sphere yields unequal oriented periods with a common radial coefficient.
\end{abstract}



\section{Introduction}\label{sec-intro}

Curvature, Ricci curvature, and Einstein conditions for Finsler warped products are studied in \cite{KozmaPeterVarga2001,ChenShenZhao2018,WuZhangSunHan2023}, while radial conformal transformations are considered in \cite{RenZhangZhao2023}. Fermat, Randers, and Finsler optical structures associated with null geodesics of static and stationary spacetimes are developed in \cite{GibbonsHerdeiroWarnickWerner2009}. Constant-radius null orbits, optical curvature, and radial stability are examined in \cite{ClaudelVirbhadraEllis2001,CardosoEtAl2009,Qiao2022,Caponio2026,QiaoLiXieGuo2025}.

We study the diagonal conic Finsler structure
$F^{2}=b(r)(v^{1})^{2}+a(r)\breve\alpha(u,\breve v)^{2}$ on
$M=I\times\breve M$, where $a,b>0$ and $\breve\alpha$ is an arbitrary
angular Finsler norm. When $\breve\alpha$ is Riemannian, the function
$\breve\alpha^{2}$ is quadratic in $\breve v$, the metric extends
smoothly to every nonzero tangent direction, and the radial-pole and
angular-pole descriptions yield the same radial sectional curvature.
For a non-Riemannian angular norm, $\breve\alpha^{2}$ is not necessarily
twice differentiable at $\breve v=0$, so the smooth geometry is
considered on the cone
$\Cone=\{(r,u;v^1,\breve v)\in TM:\breve v\neq0\}$. The angular-pole
radial-transverse flag remains defined on this domain and provides the
corresponding Finsler description of the radial Jacobi mode. We obtain
the spray curvature endomorphism and its Ricci contraction. At a
critical radius, we identify the angular-pole radial-transverse flag
curvature with the negative optical-arclength radial Jacobi coefficient;
its sign distinguishes exponential and oscillatory radial behavior,
while the zero coefficient gives affine radial solutions. We prove that
radial conformal factors preserve the diagonal class and that a
compactly supported conformal two-jet determines a selected critical
radius and an arbitrary signed radial coefficient. To the best of our
knowledge, the flag-Jacobi identity for an arbitrary angular Finsler
norm and the compactly supported two-jet assignment are new for this
diagonal conic class. The navigation Randers sphere provides two
oriented closed angular geodesics with unequal periods and the same
assigned radial coefficient.
\section{Preliminaries}
\label{sec-prelim}

We follow the metric, spray, curvature, and flag-curvature conventions in \cite{BaoChernShen2000} and the conic convention,\cite{JavaloyesSanchez2014}. Capital indices $A,B,C$ range over $1,\ldots,n$, lower-case indices $a,b,c$ range over $2,\ldots,n$, and repeated indices are summed. Every fiber derivative is taken on
\begin{equation*}
 \Cone:=\{(r,u;v^1,\breve v)\in TM:\breve v\neq0\}.
\end{equation*}
All radial functions are smooth on the interval under consideration.

For a Finsler function $F$, the fundamental tensor and canonical spray coefficients are
\[
 g_{AB}=\frac12(F^2)_{y^Ay^B},
 \qquad
 G^A=\frac14g^{AB}\bigl((F^2)_{x^Cy^B}y^C-(F^2)_{x^B}\bigr).
\]
Euler homogeneity yields $(F^2)_{y^A}y^A=2F^2$ and $g^{AB}(F^2)_{y^B}=2y^A$. The spray curvature endomorphism is
\begin{equation}\label{eq-Riem-spray}
 R^A{}_B
 =2\partial_{x^B}G^A-y^C\partial_{x^Cy^B}G^A
 +2G^C\partial_{y^Cy^B}G^A
 -\partial_{y^C}G^A\partial_{y^B}G^C.
\end{equation}
Its Ricci trace is
\[
 \Ric_F=R^A{}_A.
\]
For a flag with pole $y$ and transverse vector $V$, the flag curvature is
\begin{equation}\label{eq-flag-def}
 \mathbf K(x,y;V)
 =\frac{g_{AB}R^A{}_C(x,y)V^BV^C}
 {g_y(y,y)g_y(V,V)-g_y(y,V)^2}.
\end{equation}

We write $M=I\times\breve M$, with coordinates $(r,u^2,\ldots,u^n)$ and fiber coordinates $(v^1,\breve v^2,\ldots,\breve v^n)$. The angular factor $(\breve M,\breve\alpha)$ has fundamental tensor $\breve g_{ab}$, spray $\breve G^a$, spray curvature endomorphism $\breve R^a{}_b$, and Ricci scalar $\breve\Ric$. Product Finsler structures and the corresponding curvature notation are studied in \cite{KozmaPeterVarga2001,ChenShenZhao2018,WuZhangSunHan2023}.

\section{Curvature, Radial Jacobi Fields, and Conformal Two-Jet Control}
\label{sec-curv}

We study the diagonal conic structure on $\Cone$. We obtain its conic regularity, fundamental tensor, spray, curvature block, Ricci contraction, and distinguished flag curvatures. We examine radial conformal scaling, constant-radius angular geodesics, radial Jacobi fields, signed radial coefficients, local two-jet control, and the nonreversible Randers example. We also compare the conic formulas with the corresponding Riemannian optical identities.

\begin{definition}\label{def-Wq}
Let $M=I\times\breve M$, let $a,b\in C^\infty(I)$ be positive, and let $\breve\alpha$ be a Finsler norm on $\breve M$. Write
\[
 \Cone=\{(r,u;v^1,\breve v)\in TM:\breve v\neq0\}.
\]
The function
\begin{equation}\label{eq-Wq-F2}
 F^2=b(r)(v^1)^2+a(r)\breve\alpha(u,\breve v)^2
\end{equation}
on $\Cone$ is called a \emph{diagonal conic warped-product Finsler structure}. It is in \emph{normalized form} when $b\equiv1$, in which case
\[
 F^2=(v^1)^2+h(r)\breve\alpha(u,\breve v)^2,
 \qquad h=a>0.
\]
If $\breve\alpha$ is Riemannian, the same formula extends to a regular Finsler metric on $TM\setminus\{0\}$.
\end{definition}

The conic regularity convention is described for conic Finsler metrics, \cite{JavaloyesSanchez2014}. A positively two-homogeneous function on a finite-dimensional vector space with a $C^2$ extension at the origin is a quadratic form, \cite[Remark~5(1)]{BernalJavaloyesSanchez2020}.

\begin{proposition}\label{prop-conic-domain}
Let
\[
 F^2=b(r)(v^1)^2+a(r)\breve\alpha(u,\breve v)^2,
 \qquad a(r)>0,\qquad b(r)>0,
\]
be defined on
\[
 \Cone=\{(r,u;v^1,\breve v)\in TM:\breve v\neq0\}.
\]
Then $F$ is smooth and strongly convex on $\Cone$. If $\breve\alpha$ is Riemannian, $F$ extends to a regular Finsler metric on $TM\setminus\{0\}$. Conversely, if $F^2$ is twice fiber differentiable at every vector $(v^1,0)$ with $v^1\neq0$, then $\breve\alpha(u,\breve v)^2$ is quadratic in $\breve v$ for every $u\in\breve M$, and hence $\breve\alpha$ is Riemannian.
\end{proposition}

\begin{proof}
On $\Cone$, the angular norm is evaluated only at $\breve v\neq0$, where $\breve\alpha$ is smooth and strongly convex.  The fundamental tensor is block diagonal with entries $b(r)$ and $a(r)\breve g_{ab}$, so it is positive definite because $a,b>0$ and $\breve g_{ab}$ is positive definite.  This proves smoothness and strong convexity on $\Cone$.

If $\breve\alpha$ is Riemannian, then $\breve\alpha^2$ is a quadratic form in $\breve v$, and \eqref{eq-Wq-F2} is a smooth positive quadratic form in every nonzero fiber direction; hence it defines a regular Finsler metric on $TM\setminus\{0\}$.

Conversely, we fix $r\in I$, $u\in\breve M$, and $v^1\neq0$.  Twice fiber differentiability of $F^2$ at the pure radial vector $(v^1,0)$ implies that the angular variation
\[
 H_u(w):=\frac{1}{a(r)}\left(F^2(r,u;v^1,w)-b(r)(v^1)^2\right)
\]
extends as a $C^2$ function at $w=0$.  On $w\neq0$ this function is $\breve\alpha(u,w)^2$ and is positively two-homogeneous.  The two-homogeneous regularity result gives that $\breve\alpha(u,w)^2$ is quadratic in $w$.  Since $u$ was arbitrary, the angular norm is Riemannian.
\end{proof}

The class in \eqref{eq-Wq-F2} is the specialization
\begin{equation*}
 F=\breve\alpha\,\phi(r,s),
 \qquad
 s=\frac{v^1}{\breve\alpha},
 \qquad
 \phi(r,s)=\sqrt{a(r)+b(r)s^2},
\end{equation*}
of the general Finsler warped products studied in
\cite{ChenShenZhao2018}. Its elementary derivatives are
\begin{align}
 \phi_r&=\frac{a'+b's^2}{2\phi},
 &
 \phi_s&=\frac{bs}{\phi},
 &
 \phi_{ss}&=\frac{ab}{\phi^3},\label{eq-phi-derivatives}\\
 \phi_{rs}&=
 \frac{s\bigl(2ab'-a'b+bb's^2\bigr)}{2\phi^3}.
\end{align}
In the standard notation
\begin{equation*}
 \Psi=
 \frac{\phi_r(\phi_s+s\phi_{ss})-s\phi_s\phi_{rs}}
 {2\phi\phi_{ss}},
 \qquad
 \Theta=-\frac{\phi_r-s\phi_{rs}}{2\phi_{ss}},
\end{equation*}
substitution of \eqref{eq-phi-derivatives} gives
\begin{equation}\label{eq-Psi-Theta-specialized}
 \Psi=\frac{a'}{2a}s,
 \qquad
 \Theta=-\frac{a'}{4b}
 +\left(\frac{b'}{4b}-\frac{a'}{2a}\right)s^2.
\end{equation}
Thus the standard spray form
\[
 G^1=\breve\alpha^2(\Theta+s\Psi),
 \qquad
 G^a=\breve G^a+\breve\alpha\Psi\,\breve v^a
\]
reduces to Lemma~\ref{lem-spray-Wq}. After this specialization, the radial combinations entering the expanded curvature endomorphism are
\begin{equation}\label{eq-literature-PQ-map}
 P=\frac{a''}{2b}-\frac{a'b'}{4b^2}-\frac{(a')^2}{4ab},
 \qquad
 Q=\frac{(a')^2}{4ab}.
\end{equation}
The formulas \eqref{eq-Psi-Theta-specialized} and \eqref{eq-literature-PQ-map} specialize the warped-product quantities, \cite{ChenShenZhao2018,WuZhangSunHan2023}. Direct calculation on $\Cone$ retains the intrinsic angular term $\breve R^{a}{}_{b}$ for an arbitrary angular Finsler norm $\breve\alpha$.

\begin{lemma}\label{lem-block}
Let
\[
 F^2=b(r)(v^1)^2+a(r)\breve\alpha(u,\breve v)^2,
 \qquad a(r)>0,\qquad b(r)>0,
\]
be defined on $\Cone$. Its fundamental tensor satisfies
\[
 g_{11}=b(r),\qquad g_{1a}=0,\qquad
 g_{ab}=a(r)\breve g_{ab}(u,\breve v),
\]
and its inverse tensor satisfies
\[
 g^{11}=\frac1{b(r)},\qquad g^{1a}=0,\qquad
 g^{ab}=\frac1{a(r)}\breve g^{ab}(u,\breve v).
\]
\end{lemma}

\begin{proof}
The function $b(r)(v^1)^2$ depends on $v^1$ and the function $a(r)\breve\alpha^2$ depends on the angular fiber variables. Direct differentiation yields
\begin{align*}
 g_{11}
 &=\frac12\frac{\partial^2}{\partial (v^1)^2}\left(b(r)(v^1)^2+a(r)\breve\alpha^2\right)=b(r),\\
 g_{1a}
 &=\frac12\frac{\partial^2}{\partial v^1\partial \breve v^a}\left(b(r)(v^1)^2+a(r)\breve\alpha^2\right)=0,\\
 g_{ab}
 &=\frac12\frac{\partial^2}{\partial \breve v^a\partial \breve v^b}\left(b(r)(v^1)^2+a(r)\breve\alpha^2\right)=a(r)\breve g_{ab}.
\end{align*}
Thus the matrix of $g_{AB}$ is block diagonal. Inverting the scalar radial block and the angular block yields $g^{11}=1/b$, $g^{1a}=0$, and $g^{ab}=\breve g^{ab}/a$.
\end{proof}

\begin{lemma}\label{lem-spray-Wq}
Let
\[
 F^2=b(r)(v^1)^2+a(r)\breve\alpha(u,\breve v)^2,
 \qquad a(r)>0,\qquad b(r)>0,
\]
be defined on $\Cone$, where $\breve G^a$ are the spray coefficients of $\breve\alpha$. Then the spray coefficients of $F$ are
\begin{equation}\label{eq-spray-Wq}
 G^1=\frac{b'(r)}{4b(r)}(v^1)^2-\frac{a'(r)}{4b(r)}\breve\alpha^2,
 \qquad
 G^a=\breve G^a+\frac{a'(r)}{2a(r)}v^1\breve v^a,
 \quad a\in\{2,\ldots,n\}.
\end{equation}
\end{lemma}

\begin{proof}
We write $\mathcal Q(u,\breve v)=\breve\alpha(u,\breve v)^2$. From Lemma~\ref{lem-block}, we have
\[
 g^{11}=\frac1b,\qquad g^{1a}=0,\qquad g^{ab}=\frac1a\breve g^{ab}.
\]
Only the coordinate $x^1=r$ occurs in $a$ and $b$, while $\mathcal Q$ depends only on
$(u,\breve v)$.  Consequently,
\[
 (F^2)_{r}=b'(v^1)^2+a'\mathcal Q,
 \qquad
 (F^2)_{u^j}=a\mathcal Q_{u^j},
\]
and
\[
 (F^2)_{r v^1}=2b'v^1,
 \qquad
 (F^2)_{u^j v^1}=0,
\]
\[
 (F^2)_{r\breve v^j}=a'\mathcal Q_j,
 \qquad
 (F^2)_{u^c\breve v^j}=a\mathcal Q_{u^c\breve v^j},
 \qquad \mathcal Q_j:=(\mathcal Q)_{\breve v^j}.
\]
For the radial component, the block form of $g^{AB}$ yields
\begin{align*}
 4G^1
 &=g^{11}\left((F^2)_{x^C v^1}y^C-(F^2)_r\right)\\
 &=\frac1b\left(2b'(v^1)^2-b'(v^1)^2-a'\mathcal Q\right),
\end{align*}
which is the radial identity in \eqref{eq-spray-Wq}.

For an angular component, the inverse tensor contributes through its angular entries:
\begin{align*}
 4G^a
 &=\frac1a\breve g^{aj}
 \left(a\mathcal Q_{u^c\breve v^j}\breve v^c+a'v^1\mathcal Q_j-a\mathcal Q_{u^j}\right)\\
 &=\breve g^{aj}
 \left(\mathcal Q_{u^c\breve v^j}\breve v^c-\mathcal Q_{u^j}\right)
 +\frac{a'}a v^1\breve g^{aj}\mathcal Q_j.
\end{align*}
The angular spray satisfies
\[
 4\breve G^a
 =\breve g^{aj}
 \left(\mathcal Q_{u^c\breve v^j}\breve v^c-\mathcal Q_{u^j}\right).
\]
Euler homogeneity for $\mathcal Q$ yields $\mathcal Q_j=2\breve g_{jk}\breve v^k$, and hence
$\breve g^{aj}\mathcal Q_j=2\breve v^a$.  Therefore
\[
 4G^a=4\breve G^a+2\frac{a'}a v^1\breve v^a,
\]
which proves the angular identity in \eqref{eq-spray-Wq}.
\end{proof}

\begin{proposition}\label{prop-non-riem-angular-sector}
Let
\[
 F^2=b(r)(v^1)^2+a(r)\breve\alpha(u,\breve v)^2,
 \qquad a(r)>0,\qquad b(r)>0,
\]
be defined on $\Cone$. Let
\[
 C^F_{ABC}=\frac12\frac{\partial g_{AB}}{\partial y^C}
\]
be the Cartan tensor of $F$, and let $\breve C_{abc}$ be the Cartan tensor of $\breve\alpha$. Then
\[
 C^F_{111}=C^F_{11a}=C^F_{1ab}=0,
 \qquad
 C^F_{abc}=a(r)\breve C_{abc}.
\]
Consequently, $F$ is Riemannian if and only if $\breve\alpha$ is Riemannian. On every open region where $\breve C_{abc}\neq0$, the structure differs from a Riemannian warped product with the same radial-angular splitting and angular data.
\end{proposition}

\begin{proof}
From Lemma~\ref{lem-block}, we have $g_{11}=b(r)$, $g_{1a}=0$, and
$g_{ab}=a(r)\breve g_{ab}(u,\breve v)$.  Differentiating these identities with
respect to the fiber variables shows that every Cartan component containing a radial
index is zero and yields
\[
 C^F_{abc}
 =\frac12\frac{\partial}{\partial \breve v^c}
 \bigl(a(r)\breve g_{ab}(u,\breve v)\bigr)
 =a(r)\breve C_{abc}.
\]
Thus the Cartan tensor of $F$ vanishes only when the angular Cartan tensor
vanishes.  By the standard Finsler criterion, vanishing Cartan tensor is
equivalent to the square of the norm being quadratic in the fiber variable,
\cite{BaoChernShen2000}.  Therefore a nonzero angular Cartan tensor gives a
non-Riemannian conic optical structure.  A Riemannian warped product has an
angular metric tensor independent of angular direction; hence it differs from $a(r)\breve g_{ab}(u,\breve v)$ on every open region where
$\breve C_{abc}\neq0$.
\end{proof}

\subsection{Spray Curvature Endomorphism and Ricci Splitting}

We compute the four curvature blocks from the spray and obtain the Ricci contraction. The radial functions $P$ and $Q$ isolate the coefficients that also occur in the distinguished flag curvatures.

\begin{theorem}\label{thm-curvature-block}
Let
\[
 F^2=b(r)(v^1)^2+a(r)\breve\alpha(u,\breve v)^2,
 \qquad a(r)>0,\qquad b(r)>0,
\]
be a diagonal conic warped-product Finsler structure on $\Cone\subset TM$, where $\dim M=n$. Let $\breve R^a{}_b(u,\breve v)$ be the spray curvature endomorphism of $(\breve M,\breve\alpha)$, and write $\breve v_b=\breve g_{bc}\breve v^c$. Define
\[
 P(r)=\frac{a''(r)}{2b(r)}-\frac{a'(r)b'(r)}{4b(r)^2}
 -\frac{a'(r)^2}{4a(r)b(r)},
 \qquad
 Q(r)=\frac{a'(r)^2}{4a(r)b(r)}.
\]
Then the spray curvature endomorphism $R^A{}_B$ has the block form
\begin{align*}
 R^1{}_1&=-P(r)\breve\alpha^2,\\
 R^1{}_a&=P(r)v^1\breve v_a,\\
 R^a{}_1&=\frac{b(r)}{a(r)}P(r)v^1\breve v^a,\\
 R^a{}_b&=\breve R^a{}_b
 -\left(\frac{b(r)}{a(r)}P(r)(v^1)^2+Q(r)\breve\alpha^2\right)\delta^a_b
 +Q(r)\breve v^a\breve v_b.
\end{align*}
Consequently,
\begin{equation}\label{eq-Ricci-decomp}
 \Ric_F(x,y)=\breve\Ric(u,\breve v)
 -(n-1)\frac{b(r)}{a(r)}P(r)(v^1)^2
 -\bigl(P(r)+(n-2)Q(r)\bigr)\breve\alpha^2,
\end{equation}
where $\breve\Ric(u,\breve v)=\breve R^a{}_a(u,\breve v)$.
\end{theorem}

\begin{proof}
We write
\[
 A(u,\breve v):=\breve\alpha(u,\breve v)^2,
 \qquad
 p:=\frac{b'}{4b},
 \qquad
 q:=-\frac{a'}{4b},
 \qquad
 c:=\frac{a'}{2a}.
\]
We rewrite Lemma~\ref{lem-spray-Wq} as
\begin{equation}\label{eq-curv-proof-spray}
 G^1=p(v^1)^2+qA,
 \qquad
 G^a=\breve G^a+cv^1\breve v^a.
\end{equation}
We write $A_a=A_{\breve v^a}$ and $A_{ab}=A_{\breve v^a\breve v^b}$.
Euler homogeneity yields
\begin{equation}\label{eq-curv-proof-euler}
 A_a\breve v^a=2A,
 \qquad
 A_{ab}\breve v^a=A_b,
 \qquad
 A_a=2\breve g_{ab}\breve v^b=2\breve v_a.
\end{equation}
The two-homogeneity of $\breve G^a$ yields
\begin{equation}\label{eq-curv-proof-Ghom}
 \breve G^a_{\breve v^b}\breve v^b=2\breve G^a,
 \qquad
 \breve v^c\breve G^a_{\breve v^c\breve v^b}
 =\breve G^a_{\breve v^b}.
\end{equation}
The angular spray formula and $A_{ab}=2\breve g_{ab}$ also give
\begin{equation}\label{eq-curv-proof-spray-id1}
 2\breve G^cA_{cb}
 =\breve v^dA_{u^d\breve v^b}-A_{u^b}.
\end{equation}
We also contract the angular spray formula with $A_c=2\breve g_{ce}\breve v^e$.
Euler homogeneity of $A_{u^d}$ yields
\[
 2A_c\breve G^c=A_{u^d}\breve v^d.
\]
Differentiation in $\breve v^b$, followed by \eqref{eq-curv-proof-spray-id1}, gives
\begin{equation}\label{eq-curv-proof-spray-id2}
 A_c\breve G^c_{\breve v^b}=A_{u^b}.
\end{equation}

The derivatives of the radial spray component in \eqref{eq-curv-proof-spray} are
\[\begin{gathered}
 G^1_r=p'(v^1)^2+q'A,\quad G^1_{u^b}=qA_{u^b},\quad G^1_{v^1}=2pv^1,\quad G^1_{\breve v^b}=qA_b,\\
 G^1_{rv^1}=2p'v^1,\quad G^1_{r\breve v^b}=q'A_b,\quad G^1_{u^c\breve v^b}=qA_{u^c\breve v^b},\quad G^1_{v^1v^1}=2p,\quad G^1_{\breve v^c\breve v^b}=qA_{cb}.
\end{gathered}\]
The derivatives of the angular component are
\[\begin{gathered}
 G^a_r=c'v^1\breve v^a,\quad G^a_{u^b}=\breve G^a_{u^b},\quad G^a_{v^1}=c\breve v^a,\quad G^a_{\breve v^b}=\breve G^a_{\breve v^b}+cv^1\delta^a_b,\\
 G^a_{rv^1}=c'\breve v^a,\quad G^a_{r\breve v^b}=c'v^1\delta^a_b,\quad G^a_{v^1\breve v^b}=c\delta^a_b,\quad G^a_{\breve v^c\breve v^b}=\breve G^a_{\breve v^c\breve v^b}.
\end{gathered}\]

Substitution into \eqref{eq-Riem-spray} produces the four blocks.  For the radial-radial
entry,
\begin{align*}
 R^1{}_1
 &=2\bigl(p'(v^1)^2+q'A\bigr)-2p'(v^1)^2
   +4p\bigl(p(v^1)^2+qA\bigr)\\
 &\quad -(2pv^1)^2-qcA_a\breve v^a\\
 &=2\bigl(q'+2pq-qc\bigr)A,
\end{align*}
where \eqref{eq-curv-proof-euler} was applied in the final line.

For the radial-angular entry,
\begin{align*}
 R^1{}_b
 &=2qA_{u^b}-q'v^1A_b-q\breve v^cA_{u^c\breve v^b}
   +2q\breve G^cA_{cb}+2qcv^1\breve v^cA_{cb}\\
 &\quad -2pqv^1A_b-qA_c\breve G^c_{\breve v^b}-qcv^1A_b.
\end{align*}
The terms containing angular derivatives cancel by
\eqref{eq-curv-proof-spray-id1} and \eqref{eq-curv-proof-spray-id2}.
The identity $\breve v^cA_{cb}=A_b$ then yields
\begin{equation}\label{eq-curv-proof-R1b}
 R^1{}_b=\bigl(-q'+qc-2pq\bigr)v^1A_b.
\end{equation}

For the angular-radial entry,
\begin{align*}
 R^a{}_1
 &=2c'v^1\breve v^a-c'v^1\breve v^a
   +2c\bigl(\breve G^a+cv^1\breve v^a\bigr)\\
 &\quad -2pcv^1\breve v^a
   -c\breve G^a_{\breve v^d}\breve v^d-c^2v^1\breve v^a\\
 &=\bigl(c'+c^2-2pc\bigr)v^1\breve v^a,
\end{align*}
where \eqref{eq-curv-proof-Ghom} cancels the two terms containing
$\breve G^a$.

For the angular-angular entry, direct substitution yields
\begin{align*}
 R^a{}_b
 &=2\breve G^a_{u^b}-\breve v^d\breve G^a_{u^d\breve v^b}
   +2\breve G^d\breve G^a_{\breve v^d\breve v^b}
   -\breve G^a_{\breve v^d}\breve G^d_{\breve v^b}\\
 &\quad +\bigl(-c'+2pc-c^2\bigr)(v^1)^2\delta^a_b
   +2cqA\delta^a_b-cq\breve v^aA_b.
\end{align*}
The leading terms give $\breve R^a{}_b$ by the angular version of
\eqref{eq-Riem-spray}.  Therefore
\begin{equation}\label{eq-curv-proof-Rab}
 R^a{}_b
 =\breve R^a{}_b
 +\bigl(-c'+2pc-c^2\bigr)(v^1)^2\delta^a_b
 +2cqA\delta^a_b-cq\breve v^aA_b.
\end{equation}

The radial coefficients now simplify as follows:
\begin{align*}
 2(q'+2pq-qc)
 &=2\left(-\frac{a''}{4b}+\frac{a'b'}{4b^2}
 -\frac{a'b'}{8b^2}+\frac{(a')^2}{8ab}\right)=-P,\\
 -q'+qc-2pq
 &=\frac{a''}{4b}-\frac{a'b'}{8b^2}
 -\frac{(a')^2}{8ab}=\frac{P}{2},\\
 c'+c^2-2pc
 &=\frac{a''}{2a}-\frac{(a')^2}{4a^2}
 -\frac{a'b'}{4ab}=\frac{b}{a}P,\\
 -c'+2pc-c^2&=-\frac{b}{a}P,\\
 2cq&=-\frac{(a')^2}{4ab}=-Q.
\end{align*}
Since $A_b=2\breve v_b$, equations
\eqref{eq-curv-proof-R1b} and \eqref{eq-curv-proof-Rab} become the stated block formulas.

Finally,
\[
 \breve v^a\breve v_a=A,
\]
so tracing the angular block yields
\[
 R^a{}_a
 =\breve\Ric-(n-1)\frac{b}{a}P(v^1)^2-(n-2)QA.
\]
Together with $R^1{}_1=-PA$, this identity proves \eqref{eq-Ricci-decomp}.
\end{proof}

The supporting calculations verify the spray-kernel identity, metric self-adjointness, the Riemannian reduction, and the Einstein specialization of the curvature block.

\subsection{Distinguished Flag Curvatures}

We compute the angular flag, the angular-pole radial-transverse flag, and the radial-pole flag in the Riemannian angular case. The angular-pole radial-transverse flag lies in $\Cone$ for every angular Finsler norm and provides the curvature coefficient in the radial Jacobi equation.

\begin{theorem}\label{thm-flag}
Let
\[
 F^2=b(r)(v^1)^2+a(r)\breve\alpha(u,\breve v)^2,
 \qquad a(r)>0,\qquad b(r)>0,
\]
be defined on $\Cone$, where $\dim M=n$, and write
\[
 P(r)=\frac{a''(r)}{2b(r)}-\frac{a'(r)b'(r)}{4b(r)^2}
 -\frac{a'(r)^2}{4a(r)b(r)}.
\]
Then the following flag-curvature identities hold.
\begin{enumerate}[label=\emph{(\roman*)}]
\item Assume $n\ge3$. Let $y=(0,\breve v)$ and $V=(0,\breve W)$, where $\breve v$ and $\breve W$ are linearly independent, and let $\breve{\mathbf K}(u,\breve v;\breve W)$ be the angular flag curvature. Then
\begin{equation}\label{eq-general-angular-flag}
 \mathbf K_F\bigl((r,u);y;V\bigr)
 =\frac{1}{a(r)}\breve{\mathbf K}(u,\breve v;\breve W)
 -\frac{a'(r)^2}{4a(r)^2b(r)}.
\end{equation}
\item For every nonzero $\breve v$, the flag with $y=(0,\breve v)$ and $V=(1,0)$ satisfies
\begin{equation}\label{eq-general-radial-transverse}
 \mathbf K_F\bigl((r,u);(0,\breve v);(1,0)\bigr)
 =-\frac{P(r)}{a(r)}.
\end{equation}
This identity holds for every strongly convex angular Finsler norm.
\item If $\breve\alpha$ is Riemannian, then for $y=(v^1,0)$ with $v^1\neq0$ and $V=(0,\breve W)$ with $\breve W\neq0$,
\begin{equation}\label{eq-general-radial-flag}
 \mathbf K_F\bigl((r,u);(v^1,0);(0,\breve W)\bigr)
 =-\frac{P(r)}{a(r)}.
\end{equation}
\end{enumerate}
If $a'(r_0)=0$, then
\begin{equation}\label{eq-general-radial-ring}
 \mathbf K_F\bigl((r_0,u);(0,\breve v);(1,0)\bigr)
 =-\frac{a''(r_0)}{2a(r_0)b(r_0)}.
\end{equation}
Under the Riemannian assumption, the radial-pole flag in \emph{(iii)} has the same value.
\end{theorem}

\begin{proof}
For $y=(0,\breve v)$, we obtain from Theorem~\ref{thm-curvature-block} that
\[
 R^a{}_b
 =\breve R^a{}_b-Q\breve\alpha^2\delta^a_b
 +Q\breve v^a\breve v_b.
\]
We write
\[
 \alpha_0^2:=\breve\alpha^2,
 \qquad
 W^2:=\breve g_{\breve v}(\breve W,\breve W),
 \qquad
 c:=\breve g_{\breve v}(\breve v,\breve W),
 \qquad
 D:=\alpha_0^2W^2-c^2.
\]
Linear independence and strong convexity give $D>0$. By the definition of the angular flag curvature,
\[
 \breve g_{ac}\breve W^a\breve R^c{}_b\breve W^b
 =\breve{\mathbf K}D.
\]
Since $\breve v_b\breve W^b=c$,
\[
 g_{ac}\breve W^aR^c{}_b\breve W^b
 =a\left(\breve{\mathbf K}D-Q\alpha_0^2W^2+Qc^2\right)=a\left(\breve{\mathbf K}-Q\right)D,
 \qquad
 g_y(y,y)g_y(V,V)-g_y(y,V)^2=a^2D.
\]
By \eqref{eq-flag-def},
\[
 \mathbf K_F=\frac{a\left(\breve{\mathbf K}-Q\right)D}{a^2D}
 =\frac{\breve{\mathbf K}}a-\frac Qa
 =\frac{\breve{\mathbf K}}a-\frac{(a')^2}{4a^2b},
\]
which gives \eqref{eq-general-angular-flag} without imposing orthogonality on the chosen transverse representative.

For the angular-pole, radial-transverse flag,
\[
 g_y(y,y)=a\breve\alpha^2,\qquad g_y(V,V)=b,\qquad g_y(y,V)=0,\qquad R^1{}_1=-P\breve\alpha^2,
\]
so \eqref{eq-flag-def} gives
\[
 \mathbf K_F=\frac{-bP\breve\alpha^2}{ab\breve\alpha^2}=-\frac Pa,
\]
which is \eqref{eq-general-radial-transverse}.

We assume that the angular norm is Riemannian and choose $y=(v^1,0)$. The angular curvature
term and every term containing $\breve v$ vanish at this pole, so
\[
 R^a{}_b=-\frac baP(v^1)^2\delta^a_b.
\]
For $V=(0,\breve W)$, \eqref{eq-flag-def} gives
\[
 \mathbf K_F=\frac{a\breve g_{ac}\breve W^a\left(-\frac baP(v^1)^2\delta^c_b\right)\breve W^b}{b(v^1)^2\,aW^2}
 =\frac{-bP(v^1)^2W^2}{ab(v^1)^2W^2}=-\frac Pa,
\]
which gives \eqref{eq-general-radial-flag}. At a critical radius,
$a'(r_0)=0$, and the definition of $P$ gives
\[
 P(r_0)=\frac{a''(r_0)}{2b(r_0)}.
\]
Substitution yields \eqref{eq-general-radial-ring} for the angular-pole flag and, under the stated Riemannian assumption, for the radial-pole flag.
\end{proof}

\begin{remark}\label{rem-angular-pole-replacement}
For a non-Riemannian angular norm, the smooth Finsler calculations are confined to $\breve v\neq0$. The angular-pole radial-transverse flag in Theorem~\ref{thm-flag}\emph{(ii)} lies in $\Cone$ and represents the radial Jacobi mode. The same flag enters the Jacobi identity and the local conformal two-jet construction.
\end{remark}

\subsection{Radial Conformal Scaling}
\label{sec-rconf}

We examine radial conformal factors that preserve the product splitting and the angular Finsler norm. The transformed coefficients and the invariant ratio $b/a$ follow from the diagonal form. The conformal spray identity, \cite{BaoChernShen2000,RenZhangZhao2023} for $\bar F=e^{\sigma(x)}F$ is
\begin{equation*}
 \bar{G}^{A}=G^{A}+\sigma_0 y^{A}-\tfrac{1}{2}F^{2}\sigma^{A},
 \qquad \sigma_0:=\sigma_{x^A}y^A,\quad \sigma^A:=g^{AB}\sigma_{x^B}.
\end{equation*}
For $\sigma=\sigma(r)$,
\begin{equation*}
 \bar{G}^{A}=G^{A}+\sigma'(r)\Bigl(v^1 y^A-\tfrac{1}{2}F^2 g^{A1}\Bigr).
\end{equation*}
For a radial conformal factor, this identity preserves the diagonal form.

\begin{proposition}\label{prop-radial-conformal-preservation}
Let
\[
 F^2=b(r)(v^1)^2+a(r)\breve\alpha(u,\breve v)^2,
 \qquad a(r)>0,\qquad b(r)>0,
\]
be defined on $\Cone$ relative to a fixed product splitting and a fixed angular norm $\breve\alpha$.
\begin{enumerate}[label=\textup{(\arabic*)}]
\item For every $\sigma\in C^\infty(I)$, the function $\bar F=e^{\sigma(r)}F$ has the diagonal conic form
\[
 \bar F^2=\bar b(r)(v^1)^2+\bar a(r)\breve\alpha(u,\breve v)^2,
 \qquad
 \bar a=e^{2\sigma}a,
 \qquad
 \bar b=e^{2\sigma}b.
\]
\item Let $\rho\in C^\infty(M)$. Assume that $e^{\rho(x)}F$ has a diagonal conic form relative to the same splitting and the same angular norm, with both coefficients depending only on $r$. Then $\rho=\rho(r)$.
\end{enumerate}
\end{proposition}

\begin{proof}
For a radial function $\sigma$,
\[
 \bar F^2=e^{2\sigma(r)}F^2
 =e^{2\sigma(r)}b(r)(v^1)^2+e^{2\sigma(r)}a(r)\breve\alpha^2.
\]
Both coefficients are smooth positive functions of $r$, so the deformation is again a diagonal conic structure.

For the converse, assume that
\[
 e^{2\rho(r,u)}F^2
 =\widetilde b(r)(v^1)^2+\widetilde a(r)\breve\alpha^2
\]
on $\Cone$.  We fix $(r,u,\breve v)$ with $\breve v\neq0$.  Both sides are quadratic
polynomials in the free variable $v^1$.  Equality for every $v^1$ implies equality of their
quadratic coefficients:
\[
 e^{2\rho(r,u)}b(r)=\widetilde b(r).
\]
Since $b(r)>0$, this identity shows that $e^{2\rho(r,u)}$ is independent of $u$.
Consequently $\rho$ depends only on $r$.  The constant terms in the same polynomial
identity then give
\[
 e^{2\rho(r)}a(r)=\widetilde a(r),
\]
which is consistent with the same radial factor.  The converse follows.
\end{proof}

\begin{corollary}\label{cor-aniso-inv}
Let
\[
 F^2=b(r)(v^1)^2+a(r)\breve\alpha(u,\breve v)^2,
 \qquad a(r)>0,\qquad b(r)>0,
\]
and let $\bar F=e^{\sigma(r)}F$, where $\sigma\in C^\infty(I)$. If $\bar a=e^{2\sigma}a$ and $\bar b=e^{2\sigma}b$, then
\[
 \frac{\bar b(r)}{\bar a(r)}=\frac{b(r)}{a(r)}
 \qquad(r\in I).
\]
\end{corollary}

\begin{proof}
From Proposition~\ref{prop-radial-conformal-preservation}, we have $\bar a=e^{2\sigma}a$ and $\bar b=e^{2\sigma}b$. Therefore
\[
\frac{\bar b}{\bar a}=\frac{e^{2\sigma}b}{e^{2\sigma}a}=\frac{b}{a}.
\]
\end{proof}

\subsection{Radial and Angular Refractive Data}

We write the radial and angular optical coefficients as refractive indices. Their travel-time formulas and conservation identities lead to the constant-radius criterion and the radial Jacobi equation.

\begin{definition}\label{def-indices}
Let
\[
 F^2=b(r)(v^1)^2+a(r)\breve\alpha(u,\breve v)^2,
 \qquad a(r)>0,\qquad b(r)>0.
\]
The radial and angular refractive indices are
\[
 n_r(r)=\sqrt{b(r)},
 \qquad
 n_\Omega(r)=\sqrt{a(r)},
\]
and the optical anisotropy ratio is
\[
 \mathcal A(r)=\frac{b(r)}{a(r)}
 =\left(\frac{n_r(r)}{n_\Omega(r)}\right)^2.
\]
\end{definition}

Under $\bar{F}=e^{\sigma}F$: $\bar{n}_r=e^\sigma n_r$, $\bar{n}_\Omega=e^\sigma n_\Omega$, and
$\bar{\mathcal{A}}=\mathcal{A}$ (Corollary~\ref{cor-aniso-inv}).
Thus the conformal deformation scales both indices equally and leaves their ratio invariant.

\begin{proposition}\label{prop-travel-time}
Let $c(\lambda)=(r(\lambda),u(\lambda))$, and define
\[
 \mathcal T_F[c]=\int F(c(\lambda),\dot c(\lambda))\,d\lambda.
\]
Write $n_r(r)=\sqrt{b(r)}$ and $n_\Omega(r)=\sqrt{a(r)}$.
\begin{enumerate}[label=\emph{(\roman*)}]
\item If $\breve\alpha$ is Riemannian and $\dot u=0$, then
\[
 \mathcal T_F[c]=\int n_r(r)|\dd r|.
\]
For a non-Riemannian angular norm, the same formula is the continuous extension of the length integrand to $\breve v=0$.
\item If $r\equiv r_0$ and $\dot u\neq0$, then $c$ lies in $\Cone$ and
\[
 \mathcal T_F[c]=n_\Omega(r_0)L_{\breve\alpha}[u],
 \qquad
 L_{\breve\alpha}[u]=\int\breve\alpha(u,\dot u)\,d\lambda.
\]
\end{enumerate}
\end{proposition}

\begin{proof}
The formula $F(r,u;\dot r,0)=\sqrt{b(r)}|\dot r|$ gives the radial expression wherever the regular extension is available, and it also gives the stated continuous extension in the conic case. For a purely angular curve at $r=r_0$, one has $v^1=0$ and $\breve v\neq0$, so
\[
 F=\sqrt{a(r_0)}\,\breve\alpha(u,\dot u)=n_\Omega(r_0)\breve\alpha(u,\dot u).
\]
Integration yields the identity in \emph{(ii)}.
\end{proof}

The conservation identities required for the Jacobi calculation also occur, \cite{Caponio2026}. We write $A(u,\dot u)=\breve\alpha(u,\dot u)^2$.  From Lemma~\ref{lem-spray-Wq}, we obtain the geodesic equations
\begin{equation}\label{eq-geodesic-radial-full}
 \ddot r+\frac{b'}{2b}\dot r^2-\frac{a'}{2b}A=0
\end{equation}
and
\begin{equation}\label{eq-geodesic-angular-full}
 \ddot u^a+2\breve G^a(u,\dot u)+\frac{a'}a\dot r\dot u^a=0.
\end{equation}
For the autonomous Lagrangian $L=F^2/2$, the Euler\textendash Lagrange equations give
\[
 \frac{d}{d\lambda}
 \left(\dot x^A\frac{\partial L}{\partial\dot x^A}-L\right)=0.
\]
Two-homogeneity gives
$\dot x^A L_{\dot x^A}=2L$, so $F^2=2L$ is constant.  Denote this nonnegative constant by $E^2$.

The angular energy identity for $A/2$ is
\[
 A_{u^a}\dot u^a-2A_{\breve v^a}\breve G^a=0.
\]
Using \eqref{eq-geodesic-angular-full} and
$A_{\breve v^a}\dot u^a=2A$, one obtains
\begin{align*}
 \frac{dA}{d\lambda}
 &=A_{u^a}\dot u^a+A_{\breve v^a}\ddot u^a\\
 &=A_{u^a}\dot u^a-2A_{\breve v^a}\breve G^a
   -\frac{a'}a\dot r A_{\breve v^a}\dot u^a\\
 &=-2\frac{a'}a\dot r A.
\end{align*}
Consequently,
\[
 \frac{d}{d\lambda}(a^2A)
 =2aa'\dot r A+a^2\frac{dA}{d\lambda}=0.
\]
We write the nonnegative constant $a^2A$ as $J^2$.  Since
$E^2=b\dot r^2+aA$, substitution of $A=J^2/a^2$ gives
\begin{equation*}
 b(r)\dot r^2+\frac{J^2}{a(r)}=E^2.
\end{equation*}

For fixed $J\neq0$, we write $V_J=J^2/a$, consistent, \cite{Caponio2026,Qiao2022}.  Differentiation yields
\[
 V_J'=-\frac{J^2a'}{a^2},
 \qquad
 V_J''=-\frac{J^2a''}{a^2}+2\frac{J^2(a')^2}{a^3}.
\]
Since $a>0$ and $J\neq0$, a critical radius satisfies $V_J'(r_0)=0$ if and only if
$a'(r_0)=0$.  At such a radius the $a'$-dependent term in $V_J''$ vanishes, and therefore
\begin{equation*}
 V_J'(r_0)=0\iff a'(r_0)=0,
 \qquad
 V_J''(r_0)=-\frac{J^2a''(r_0)}{a(r_0)^2}.
\end{equation*}
Thus $a''(r_0)>0$ gives a strict local maximum of $V_J$.

\begin{corollary}\label{cor-phys-conformal-equivalence}
Let
\[
 F_i^2=b_i(r)(v^1)^2+a_i(r)\breve\alpha^2,
 \qquad a_i(r)>0,\qquad b_i(r)>0,
 \qquad i\in\{1,2\}.
\]
There exists $\sigma\in C^\infty(I)$ such that $F_2=e^{\sigma(r)}F_1$ if and only if
\[
 \frac{b_1(r)}{a_1(r)}=\frac{b_2(r)}{a_2(r)}
 \qquad(r\in I).
\]
In this case,
\[
 e^{2\sigma(r)}=\frac{a_2(r)}{a_1(r)}
 =\frac{b_2(r)}{b_1(r)}.
\]
\end{corollary}

\begin{proof}
If $F_2=e^{\sigma(r)}F_1$, then
\[
 a_2=e^{2\sigma}a_1,
 \qquad
 b_2=e^{2\sigma}b_1,
\]
and division gives $b_2/a_2=b_1/a_1$.

Conversely, assume $b_1/a_1=b_2/a_2$.  Positivity of all four coefficients implies
\[
 \frac{b_2}{b_1}=\frac{a_2}{a_1}>0.
\]
We define
\[
 \sigma(r)=\frac12\ln\frac{a_2(r)}{a_1(r)}.
\]
Then $e^{2\sigma}a_1=a_2$ and, by the ratio identity,
$e^{2\sigma}b_1=b_2$.  Hence $F_2^2=e^{2\sigma}F_1^2$.  Positivity of the two norms
gives $F_2=e^\sigma F_1$.
\end{proof}

\begin{corollary}\label{cor-phys-jet0}
Let
\[
 F^2=b(r)(v^1)^2+a(r)\breve\alpha(u,\breve v)^2,
 \qquad a(r)>0,\qquad b(r)>0,
\]
and let $\bar F=e^{\sigma(r)}F$. Fix $r_0\in I$, and let $\gamma$ be an oriented closed $\breve\alpha$-geodesic such that $c=(r_0,\gamma)$ is a constant-radius geodesic of $\bar F$. Then its spatial optical period is
\[
 \mathcal T_{\bar F}[c]
 =e^{\sigma(r_0)}\sqrt{a(r_0)}L_{\breve\alpha}[\gamma],
 \qquad
 L_{\breve\alpha}[\gamma]
 =\int\breve\alpha(\gamma,\dot\gamma)\,d\lambda.
\]
\end{corollary}

\begin{proof}
Along $c(\lambda)=(r_0,\gamma(\lambda))$ one has $\dot r=0$ and
$\dot\gamma\neq0$.  Therefore
\[
 \bar F(c,\dot c)
 =e^{\sigma(r_0)}\sqrt{a(r_0)}\,
   \breve\alpha(\gamma,\dot\gamma).
\]
Integration over one period of the closed angular geodesic yields
\[
 \mathcal T_{\bar F}[c]
 =e^{\sigma(r_0)}\sqrt{a(r_0)}
 \int\breve\alpha(\gamma,\dot\gamma)\,d\lambda,
\]
which is the stated identity.
\end{proof}

\begin{example}\label{ex-phys-sphere-period}
If $\breve M=\mathbb S^{n-1}$ and $\breve\alpha$ is the standard round norm, then every great circle has length $2\pi$. Hence every oriented closed constant-radius optical geodesic of $\bar F=e^{\sigma(r)}F$ satisfies
\[
\mathcal T_{\bar F}[(r_0,\gamma)]=2\pi e^{\sigma(r_0)}\sqrt{a(r_0)}.
\]
\end{example}

\subsection{Constant-Radius Angular Criterion}
\noindent
We characterize constant-radius angular geodesics through the radial geodesic equation. The criticality condition produces a constant-radius lift, angular closure produces an oriented closed optical geodesic, and the static Fermat correspondence gives its future null lift.

\begin{definition}\label{def-light-ring}
Let $\gamma$ be a nonconstant $\breve\alpha$-geodesic. An $F$-geodesic
\[
 c(\lambda)=(r_0,\gamma(\lambda))
\]
is called a \emph{constant-radius angular geodesic}. The number $r_0$ is called an \emph{angular critical radius} if every local nonconstant $\breve\alpha$-geodesic segment has a constant-radius lift at $r_0$. If $\gamma$ is oriented and closed, the lift is called an \emph{oriented closed constant-radius optical geodesic}. Assume that positive functions $A$ and $B$ satisfy
\[
 b(r)=\frac{B(r)}{A(r)},
 \qquad
 a(r)=\frac{r^2}{A(r)},
\]
and consider
\[
 \mathcal L=-A(r)(\dot t)^2+B(r)(\dot r)^2
 +r^2\breve\alpha(u,\dot u)^2.
\]
The future null lift of an oriented closed constant-radius optical geodesic is called a \emph{light ring}; its spatial projection is periodic, and its static time coordinate increases by the optical period.
\end{definition}

\begin{remark}\label{rem-oriented-light-ring}
For a nonreversible angular norm, each closed spatial optical orbit has a fixed orientation. Parameter reversal can lead to another oriented geodesic and another optical period. The navigation Randers example contains two oriented equatorial geodesics with unequal periods.
\end{remark}

E. Caponio \cite{Caponio2026} established the corresponding constant-radius condition for
static Lorentz\textendash Finsler spacetimes with a Finsler angular sector.  In the present positive
optical model, the condition also follows directly from the geodesic equations.  Let
$c(\lambda)=(r_0,\gamma(\lambda))$, where $\gamma$ is a nonconstant affinely
parameterized $\breve\alpha$-geodesic.  Then $\dot r=\ddot r=0$ and
$\ell=\breve\alpha(\gamma,\dot\gamma)^2>0$.  Equation
\eqref{eq-geodesic-radial-full} becomes
\[
 -\frac{a'(r_0)}{2b(r_0)}\ell=0.
\]
Since $b(r_0)>0$ and $\ell>0$, this identity is equivalent to
\begin{equation}\label{eq-angular-critical-known}
 a'(r_0)=0.
\end{equation}
Under this condition, equation \eqref{eq-geodesic-angular-full} reduces to
\[
 \ddot\gamma^a+2\breve G^a(\gamma,\dot\gamma)=0,
\]
so every local nonconstant angular geodesic segment has a constant-radius lift.  Conversely, every such
lift forces \eqref{eq-angular-critical-known}.  Since
$n_\Omega'=a'/(2\sqrt a)$, the same condition is stationarity of $n_\Omega$.
The lifted trajectory becomes an oriented closed constant-radius optical geodesic only when its angular projection is closed; a static light-ring interpretation additionally requires the static representative of Proposition~\ref{prop-static-representative}.  

For
\[
 ds^2=-A(r)dt^2+B(r)dr^2+r^2d\Omega^2,
 \qquad
 a(r)=\frac{r^2}{A(r)},\quad b(r)=\frac{B(r)}{A(r)},
\]
condition~\eqref{eq-angular-critical-known} becomes the established photon-sphere equation, \cite{ClaudelVirbhadraEllis2001,Qiao2022}
\begin{equation*}
 2A(r_0)=r_0A'(r_0).
\end{equation*}
This relation provides a Riemannian consistency check.

\subsection{The Radial Jacobi Mode}

\noindent
We compute the scalar Jacobi equation for radial variations of a constant-radius angular geodesic. The angular velocity remains nonzero, so the calculation is carried out on $\Cone$. The spray covariant derivative and Jacobi equation on conic Finsler metrics follow, \cite{JavaloyesSoares2015}.

\begin{theorem}\label{thm-jacobi-ring}
Let
\[
 F^2=b(r)(v^1)^2+a(r)\breve\alpha(u,\breve v)^2,
 \qquad a(r)>0,\qquad b(r)>0,
\]
be defined on $\Cone$. Assume $a'(r_0)=0$, and let $c(\lambda)=(r_0,\gamma(\lambda))$ be a constant-radius angular geodesic with
\[
 \ell=\breve\alpha(\gamma,\dot\gamma)^2>0.
\]
Assume that a geodesic variation $c_\varepsilon$ is defined for $|\varepsilon|<\varepsilon_0$, that every velocity $\dot c_\varepsilon$ lies in $\Cone$, and that its variational field is $J(\lambda)=\xi(\lambda)\partial_r$. Then
\begin{equation}\label{eq-jacobi-ring}
 \ddot\xi(\lambda)=\Lambda(r_0)\ell\,\xi(\lambda),
 \qquad
 \Lambda(r_0)=\frac{a''(r_0)}{2b(r_0)}.
\end{equation}
The radial coordinate field $\partial_r$ is parallel along $c$ for the spray covariant derivative.

Assume in addition that $\gamma$ is oriented and closed and that positive functions $A,B$ satisfy $b=B/A$ and $a=r^2/A$. For
\[
 \mathcal L=-A(r)(\dot t)^2+B(r)(\dot r)^2
 +r^2\breve\alpha(u,\dot u)^2,
\]
the future null lift has periodic spatial projection, and its static time increment over one circuit equals the optical period.
\end{theorem}

\begin{proof}
We define the nonlinear connection coefficients by
\[
 N^A{}_B:=\frac{\partial G^A}{\partial y^B},
\]
and write $\mathcal D/d\lambda$ for the corresponding covariant derivative along a geodesic. The identities in Lemma~\ref{lem-spray-Wq} imply
\[
 N^1{}_1=\frac{b'}{2b}v^1,
 \qquad
 N^a{}_1=\frac{a'}{2a}\breve v^a.
\]
Along $c$, one has $v^1=0$ and $a'(r_0)=0$. Hence
\[
 N^A{}_1\big|_c=0,
 \qquad
 \frac{\mathcal D}{d\lambda}\partial_r=0.
\]
Along the same geodesic, we obtain from Theorem~\ref{thm-curvature-block} that
\[
 R^1{}_1=-P(r_0)\ell,
 \qquad
 R^a{}_1=0.
\]
The criticality condition yields
\[
 P(r_0)=\frac{a''(r_0)}{2b(r_0)}=\Lambda(r_0).
\]
The Finsler Jacobi equation in spray form,\cite{JavaloyesSoares2015}, is
\[
 \frac{\mathcal D^2J}{d\lambda^2}+R_y(J)=0.
\]
Since $\partial_r$ is parallel,
\[
 \frac{\mathcal D^2J}{d\lambda^2}=\ddot\xi\,\partial_r,
 \qquad
 R_y(J)=-\Lambda(r_0)\ell\xi\,\partial_r.
\]
The radial component is therefore \eqref{eq-jacobi-ring}. The reference angular velocity is nonzero, and openness of $\Cone$ ensures that the stated sufficiently small variation remains in the domain of every fiber derivative used here.

A coordinate computation confirms the same coefficient. For a geodesic variation whose linear variation is radial, write
\[
 r_\epsilon=r_0+\epsilon\xi+O(\epsilon^2),
 \qquad
 u_\epsilon=\gamma+O(\epsilon^2).
\]
Then
\[
 \dot r_\epsilon=\epsilon\dot\xi+O(\epsilon^2),
 \qquad
 \ddot r_\epsilon=\epsilon\ddot\xi+O(\epsilon^2),
\]
and Taylor expansion at $r_0$ yields
\[
 a'(r_\epsilon)=\epsilon a''(r_0)\xi+O(\epsilon^2),
 \qquad
 \frac1{b(r_\epsilon)}=\frac1{b(r_0)}+O(\epsilon).
\]
The term $(b'/2b)\dot r_\epsilon^2$ is $O(\epsilon^2)$. Since the angular speed is $\ell+O(\epsilon^2)$, the radial geodesic equation becomes
\[
 \epsilon\ddot\xi
 -\epsilon\frac{a''(r_0)}{2b(r_0)}\ell\xi
 +O(\epsilon^2)=0.
\]
Its linear variation is \eqref{eq-jacobi-ring}. In the angular equation, the additional radial term is
\[
 \frac{a'(r_\epsilon)}{a(r_\epsilon)}\dot r_\epsilon\dot\gamma^a=O(\epsilon^2),
\]
so the radial field satisfies the full linearized geodesic system.

Under the additional static assumptions, the null condition $\mathcal L=0$ and $\dot t>0$ imply
\[
 \dot t=\sqrt{\frac{B(r)}{A(r)}\dot r^{\,2}+\frac{r^2}{A(r)}\breve\alpha(u,\dot u)^2}=F(c,\dot c).
\]
For an oriented closed $\gamma$, integration over one spatial circuit yields
\[
 \Delta t=\int F(c,\dot c)\,d\lambda=\mathcal T_F[c].
\]
Thus the spatial projection is periodic, and its static time increment equals the optical period.
\end{proof}

Circular-null-geodesic radial behavior and Lyapunov rates are studied in  \cite{CardosoEtAl2009}, while the intrinsic Randers-Finsler flag-curvature classification is developed in \cite{QiaoLiXieGuo2025}. We obtain the corresponding diagonal conic identity.

\begin{theorem}\label{thm-flag-jacobi}
Let
\[
 F^2=b(r)(v^1)^2+a(r)\breve\alpha(u,\breve v)^2,
 \qquad a(r)>0,\qquad b(r)>0,
\]
be defined on $\Cone$. Assume $a'(r_0)=0$, and let $c(\lambda)=(r_0,\gamma(\lambda))$ be a constant-radius angular geodesic with $\ell=\breve\alpha(\gamma,\dot\gamma)^2>0$. For the flag pole $y=(0,\dot\gamma)$ and radial transverse vector $V=(1,0)$,
\begin{equation}\label{eq-flag-jacobi-curvature}
 \mathbf K_F\bigl((r_0,\gamma);y;V\bigr)
 =-\frac{a''(r_0)}{2a(r_0)b(r_0)}.
\end{equation}
If $F(c,\dot c)=1$, define
\begin{equation}\label{eq-flag-lyapunov-identity}
 C_{\rm opt}(r_0)=\frac{a''(r_0)}{2a(r_0)b(r_0)}
 =-\mathbf K_F\bigl((r_0,\gamma);y;V\bigr).
\end{equation}
For a radial variational field $J=\xi\partial_r$,
\begin{equation}\label{eq-flag-jacobi-equation}
 \ddot\xi=C_{\rm opt}(r_0)\xi,
 \qquad
 \ddot\xi+\mathbf K_F\bigl((r_0,\gamma);y;V\bigr)\xi=0.
\end{equation}
If $C_{\rm opt}>0$, the real optical growth rate is
$\gamma_{\rm opt}=\sqrt{C_{\rm opt}}$. If $C_{\rm opt}<0$, the radial
frequency is $\sqrt{-C_{\rm opt}}$. If $C_{\rm opt}=0$, the radial
equation is degenerate and its solutions are $\xi(\lambda)=C_1+C_2\lambda$, where $C_1,C_2\in\RR$.
These statements hold for every strongly convex angular Finsler norm.
\end{theorem}

\begin{proof}
At the pole $y=(0,\dot\gamma)$ and for $V=(1,0)$, the block form of the fundamental tensor yields
\[
 g_y(y,y)=a(r_0)\ell,\qquad g_y(V,V)=b(r_0),\qquad g_y(y,V)=0.
\]
The curvature block satisfies
\[
 R^1{}_1(r_0,y)=-P(r_0)\ell,\qquad P(r_0)=\frac{a''(r_0)}{2b(r_0)},
\]
where the second identity follows from $a'(r_0)=0$. Hence
\[
 g_{11}R^1{}_1=-b(r_0)P(r_0)\ell,
\]
and
\[
 g_y(y,y)g_y(V,V)-g_y(y,V)^2=a(r_0)b(r_0)\ell.
\]
Division gives
\[
 \mathbf K_F\bigl((r_0,\gamma);y;V\bigr)=-\frac{P(r_0)}{a(r_0)}=-\frac{a''(r_0)}{2a(r_0)b(r_0)},
\]
which proves \eqref{eq-flag-jacobi-curvature}.

Optical-arclength normalization yields $1=a(r_0)\ell$ and therefore $\ell=1/a(r_0)$. The radial Jacobi equation becomes
\[
 \ddot\xi=\frac{a''(r_0)}{2b(r_0)}\ell\xi=\frac{a''(r_0)}{2a(r_0)b(r_0)}\xi=-\mathbf K_F\bigl((r_0,\gamma);y;V\bigr)\xi.
\]
We obtain \eqref{eq-flag-lyapunov-identity} and \eqref{eq-flag-jacobi-equation}. If $C_{\rm opt}>0$, the characteristic roots are $\pm\sqrt{C_{\rm opt}}$. If $C_{\rm opt}<0$, writing $\omega=\sqrt{-C_{\rm opt}}$ gives the roots $\pm i\omega$. If $C_{\rm opt}=0$, two integrations yield $\xi(\lambda)=C_1+C_2\lambda$.
\end{proof}

\begin{corollary}\label{cor-radial-linear-behavior}
Let
\[
 F^2=b(r)(v^1)^2+a(r)\breve\alpha(u,\breve v)^2,
 \qquad a(r)>0,\qquad b(r)>0,
\]
be defined on $\Cone$, and let $c(\lambda)=(r_0,\gamma(\lambda))$ be a constant-radius angular geodesic with $a'(r_0)=0$ and $\ell=\breve\alpha(\gamma,\dot\gamma)^2>0$. Write
\[
 \Lambda(r_0)=\frac{a''(r_0)}{2b(r_0)}.
\]
Then $\ddot\xi=\Lambda(r_0)\ell\xi$. For arbitrary constants $C_1,C_2\in\RR$, the solutions have the following forms.
\begin{enumerate}[label=\emph{(\roman*)}]
\item If $\Lambda(r_0)>0$, then
\[
 \xi(\lambda)=C_1e^{\sqrt{\Lambda(r_0)\ell}\lambda}
 +C_2e^{-\sqrt{\Lambda(r_0)\ell}\lambda}.
\]
Unless the coefficient of the increasing exponential vanishes, the radial perturbation grows exponentially.
\item If $\Lambda(r_0)<0$, then
\[
 \xi(\lambda)=C_1\cos(\omega\lambda)+C_2\sin(\omega\lambda),
 \qquad
 \omega=\sqrt{-\Lambda(r_0)\ell}.
\]
Every radial perturbation is bounded and oscillatory.
\item If $\Lambda(r_0)=0$, then $\xi(\lambda)=C_1+C_2\lambda$. The linear equation is degenerate, and its nonlinear radial behavior requires a separate analysis.
\end{enumerate}
Under optical-arclength normalization,
\[
 C_{\rm opt}=\frac{a''(r_0)}{2a(r_0)b(r_0)},
\]
and the angular-pole radial-transverse flag curvature is $-C_{\rm opt}$.
\end{corollary}

\begin{proof}
We write $\mu=\Lambda(r_0)\ell$.  Since $\ell>0$ and $b(r_0)>0$, the signs of
$\mu$, $\Lambda(r_0)$, and $a''(r_0)$ coincide.  If $\mu>0$, the general solution is
\[
 \xi(\lambda)=C_1e^{\sqrt\mu\lambda}+C_2e^{-\sqrt\mu\lambda},
\]
so generic perturbations have an exponentially increasing component.  If $\mu<0$, write
$\omega=\sqrt{-\mu}$; then
\[
 \xi(\lambda)=C_1\cos(\omega\lambda)+C_2\sin(\omega\lambda),
\]
which is bounded.  If $\mu=0$, then $\xi=C_1+C_2\lambda$, and the linearized equation
does not determine the nonlinear radial behavior.  This gives the three radial alternatives.  
\end{proof}

\begin{remark}\label{rem-classical-photon-sphere-radial-behavior}
For the Schwarzschild optical metric, $a'(3M)=0$ and $a''(3M)=18>0$ by Example~\ref{ex-schwarzschild}. Hence the exterior photon sphere has the exponential radial sign, consistently with the standard spherical calculation, \cite{ClaudelVirbhadraEllis2001}.
\end{remark}

\subsection{Signed Radial Coefficients and Lyapunov Rates}

We express the radial Jacobi coefficient for an affine parameter, optical arclength, and a normalized static time coordinate. Positive coefficients yield real growth rates, while negative coefficients yield radial frequencies.

\begin{theorem}\label{thm-lyapunov}
Let
\[
 F^2=b(r)(v^1)^2+a(r)\breve\alpha(u,\breve v)^2,
 \qquad a(r)>0,\qquad b(r)>0,
\]
be defined on $\Cone$, and write $n_r=\sqrt b$ and $n_\Omega=\sqrt a$. Fix $r_0\in I$ with $a'(r_0)=0$, and let $c(\lambda)=(r_0,\gamma(\lambda))$ be a constant-speed angular geodesic with $\ell=\breve\alpha(\gamma,\dot\gamma)^2>0$. Define
\begin{equation}\label{eq-lyapunov}
 C_{\rm aff}(r_0)=\frac{a''(r_0)}{2b(r_0)}\ell.
\end{equation}
The radial equation is $\ddot\xi=C_{\rm aff}\xi$. If $C_{\rm aff}>0$, its affine growth rate is
\begin{equation}\label{eq-lyapunov-indices}
 \kappa_{\rm aff}(r_0)=\sqrt{C_{\rm aff}(r_0)}
 =\frac{\sqrt{n_\Omega(r_0)n_\Omega''(r_0)}}{n_r(r_0)}\sqrt\ell.
\end{equation}
If $C_{\rm aff}<0$, the affine radial frequency is $\sqrt{-C_{\rm aff}}$; if $C_{\rm aff}=0$, the radial equation is degenerate and its solutions are affine.

Under $F(c,\dot c)=1$, define
\begin{equation}\label{eq-lyapunov-optical}
 C_{\rm opt}(r_0)=\frac{a''(r_0)}{2a(r_0)b(r_0)}.
\end{equation}
Then $\ddot\xi=C_{\rm opt}\xi$. In a static realization with normalized time, the coordinate-time growth rate is $\lambda_t=\sqrt{C_{\rm opt}}$ when $C_{\rm opt}>0$.

Let $\bar F=e^{\sigma(r)}F$ and assume
\[
 a'(r_0)+2a(r_0)\sigma'(r_0)=0.
\]
The signed affine coefficient of $\bar F$ is
\begin{equation}\label{eq-lyapunov-dressed}
 \bar C_{\rm aff}(r_0)
 =\frac{a''(r_0)+2a(r_0)\sigma''(r_0)-4a(r_0)\sigma'(r_0)^2}{2b(r_0)}\ell.
\end{equation}
If $\sigma(r)$ tends to $\sigma_\infty$ on a static end with normalized time coordinate, then
\begin{equation}\label{eq-lyapunov-optical-dressed}
 \bar C_{\rm opt}(r_0)
 =e^{-2(\sigma(r_0)-\sigma_\infty)}
 \frac{a''(r_0)+2a(r_0)\sigma''(r_0)-4a(r_0)\sigma'(r_0)^2}{2a(r_0)b(r_0)}.
\end{equation}
For $\bar C_{\rm opt}>0$, the normalized coordinate-time growth rate is $\bar\lambda_t=\sqrt{\bar C_{\rm opt}}$; for $\bar C_{\rm opt}<0$, the normalized radial frequency is $\sqrt{-\bar C_{\rm opt}}$.
\end{theorem}

\begin{proof}
If $a'(r_0)=0$, we obtain from Theorem~\ref{thm-jacobi-ring} that
\[
 \ddot\xi=\frac{a''(r_0)}{2b(r_0)}\ell\,\xi,
\]
which gives \eqref{eq-lyapunov}. Since $n_\Omega=\sqrt a$,
\[
 n_\Omega''=\frac{a''}{2\sqrt a}-\frac{(a')^2}{4a^{3/2}}.
\]
At the critical radius,
\[
 n_\Omega(r_0)n_\Omega''(r_0)=\frac{a''(r_0)}2,
 \qquad n_r(r_0)^2=b(r_0).
\]
Taking the positive square root when $C_{\rm aff}>0$ gives
\eqref{eq-lyapunov-indices}. The remaining affine sign cases follow from the
constant-coefficient equation.

A constant rescaling of the affine parameter rescales $\ell$ and
$C_{\rm aff}$ by the inverse square factor. Optical arclength determines the scale through
\[
 1=F(c,\dot c)^2=a(r_0)\ell,
 \qquad \ell=\frac1{a(r_0)}.
\]
Substitution into \eqref{eq-lyapunov} proves \eqref{eq-lyapunov-optical}. In a
static realization, Fermat's relation identifies optical length with the normalized
coordinate-time lapse along future null trajectories.

For the deformed structure,
\[
 \bar a=e^{2\sigma}a,\qquad \bar b=e^{2\sigma}b,\qquad
 \bar a''=e^{2\sigma}\left(a''+4\sigma'a'+2a\sigma''+4a(\sigma')^2\right).
\]
The criticality condition yields $a'=-2a\sigma'$ at $r_0$, hence
\begin{equation}\label{eq-lyap-proof-bara2}
 \bar a''(r_0)=e^{2\sigma(r_0)}
 \left(a''+2a\sigma''-4a(\sigma')^2\right)_{r_0}.
\end{equation}
The affine coefficient is $\bar a''\ell/(2\bar b)$, so
\eqref{eq-lyap-proof-bara2} proves \eqref{eq-lyapunov-dressed}.

For optical arclength of $\bar F$, the signed coefficient is
\[
 \frac{\bar a''}{2\bar a\bar b}
 =e^{-2\sigma(r_0)}
 \frac{a''+2a\sigma''-4a(\sigma')^2}{2ab}.
\]
Asymptotic time normalization multiplies the Fermat norm by
$e^{-\sigma_\infty}$ and therefore multiplies this coefficient by
$e^{2\sigma_\infty}$. The resulting factor is
$e^{-2(\sigma(r_0)-\sigma_\infty)}$, which gives
\eqref{eq-lyapunov-optical-dressed}.
\end{proof}

\begin{corollary}\label{cor-signed-two-jet-assignment}
Let $a,b\in C^\infty(I)$ be positive, let
\[
 F^2=b(r)(v^1)^2+a(r)\breve\alpha(u,\breve v)^2,
\]
and let $\bar F=e^{\sigma(r)}F$ with $\sigma\in C^\infty(I)$. Assume that $\sigma$ tends to $\sigma_\infty$ on a specified static end, and write
\[
 \Delta\sigma=\sigma(r_0)-\sigma_\infty,
 \qquad
 \sigma'(r_0)=-\frac{a'(r_0)}{2a(r_0)}.
\]
For every $C_*\in\RR$, the identity
\[
 e^{-2\Delta\sigma}
 \frac{a''(r_0)+2a(r_0)\sigma''(r_0)-4a(r_0)\sigma'(r_0)^2}{2a(r_0)b(r_0)}
 =C_*
\]
determines
\begin{equation}\label{eq-sigma2-normalized}
 \sigma''(r_0)
 =\frac{2a(r_0)b(r_0)e^{2\Delta\sigma}C_*-a''(r_0)
 +4a(r_0)\sigma'(r_0)^2}{2a(r_0)}.
\end{equation}
The resulting normalized optical coefficient satisfies $\bar C_{\rm opt}(r_0)=C_*$. For fixed $\Delta\sigma$ and $\sigma'(r_0)$, this coefficient is strictly increasing as a function of $\sigma''(r_0)$. Positive and negative values of $C_*$ yield exponential and oscillatory radial modes, while $C_*=0$ yields a degenerate affine radial mode.
\end{corollary}

\begin{proof}
We write $a_0=a(r_0)$ and $b_0=b(r_0)$. Equation
\eqref{eq-lyapunov-optical-dressed} with
$\bar C_{\rm opt}(r_0)=C_*$ is equivalent to
\[
 2a_0b_0e^{2\Delta\sigma}C_*
 =a''(r_0)+2a_0\sigma''(r_0)-4a_0(\sigma'(r_0))^2.
\]
Solving for $\sigma''(r_0)$ yields \eqref{eq-sigma2-normalized}. The coefficient of
$\sigma''(r_0)$ in $\bar C_{\rm opt}$ is
\[
 e^{-2\Delta\sigma}\frac1{b_0}>0,
\]
which proves strict monotonicity.
\end{proof}

\begin{proposition}\label{prop-static-representative}
Let $A(r)>0$ and $B(r)>0$ on a static region, and consider
\begin{equation}\label{eq-static-lorentz-finsler}
 \mathcal L=-A(r)(\dot t)^2+B(r)(\dot r)^2+r^2\breve\alpha(u,\dot u)^2.
\end{equation}
For future null curves with $\dot u\neq0$, the fiber Hessian of $\mathcal L$ has Lorentzian signature, and the null condition is equivalent to
\[
 \dot t=F,
 \qquad
 F^2=\frac{B(r)}{A(r)}(\dot r)^2
 +\frac{r^2}{A(r)}\breve\alpha(u,\dot u)^2.
\]
For $\dot u=0$, the formula is interpreted by continuous extension and is smooth at that direction when $\breve\alpha$ is Riemannian. If a spatial optical geodesic closes after an optical period $T>0$, its future null lift returns to the same spatial point and satisfies $\Delta t=T$.

For $\bar F=e^{\sigma(r)}F$, a static representative is
\begin{equation}\label{eq-static-deformed-representative}
 \bar{\mathcal L}=-e^{-2\sigma(r)}A(r)(\dot t)^2
 +B(r)(\dot r)^2+r^2\breve\alpha(u,\dot u)^2.
\end{equation}
If the angular factor is Riemannian and spherical, its Lorentzian coefficients satisfy
\begin{equation}\label{eq-static-coefficients}
 \bar A=e^{-2\sigma}A,
 \qquad
 \bar B=B.
\end{equation}
If $\sigma(r)$ tends to $\sigma_\infty$ on a noncompact end with normalized time coordinate, then
\begin{equation}\label{eq-static-normalized-lapse}
 \bar A_{\rm norm}=e^{-2(\sigma-\sigma_\infty)}A.
\end{equation}
\end{proposition}

\begin{proof}
On the conic domain with $\dot u\neq0$, the fiber Hessian of \eqref{eq-static-lorentz-finsler} has one negative time entry and the positive spatial block
\[
 \operatorname{diag}\bigl(B,r^2\breve g_{ab}\bigr).
\]
Thus it has Lorentzian signature. For a future null curve, $\mathcal L=0$ and $\dot t>0$, so
\[
 \dot t
 =\sqrt{\frac BA(\dot r)^2+\frac{r^2}{A}\breve\alpha^2}
 =F.
\]
Integration over one closed spatial optical circuit gives $\Delta t=\int F\,d\lambda=T>0$. This shows the distinction between a periodic spatial projection and a closed spacetime curve.
Replacing $A$ by $e^{-2\sigma}A$ multiplies the resulting Fermat norm by $e^\sigma$, which gives \eqref{eq-static-deformed-representative} for smooth conic null directions. Pure radial null directions are covered only by the continuous extension stated in the proposition. In the spherical Riemannian case, comparison with
$-\bar A\,dt^2+\bar B\,dr^2+r^2d\Omega^2$ gives
\eqref{eq-static-coefficients}. If $\sigma$ approaches $\sigma_\infty$, rescaling the static time so that the lapse approaches one replaces $e^{-2\sigma}A$ by
$e^{-2(\sigma-\sigma_\infty)}A$, proving \eqref{eq-static-normalized-lapse}.
\end{proof}

\subsection{Local Conformal Two-Jet Control}
\label{sec-twojet}

We construct a compactly supported radial conformal function whose value, first derivative, and second derivative at a selected radius determine the criticality condition and the signed radial coefficient. The value $c_0=\sigma(r_0)$ enters the local coefficient, while $\Delta\sigma=\sigma(r_0)-\sigma_\infty$ enters the coefficient normalized at a specified static end. Support disjoint from a neighborhood of that end preserves its normalization.

\begin{proposition}\label{prop-local-conformal-twojet-flexibility}
Let $I\subset\RR$ be open, and let
\[
 F^2=b(r)(v^1)^2+a(r)\breve\alpha(u,\breve v)^2,
 \qquad a(r)>0,\qquad b(r)>0,
\]
be defined on $\Cone$ over $I$. Fix $r_0\in I$, choose a compact interval $J\Subset I$ containing $r_0$ in its interior, and fix $c_0,C_*\in\RR$. Then there exists $\sigma\in C_c^\infty(J)$ such that $\bar F=e^{\sigma(r)}F$, $\bar a=e^{2\sigma}a$, and $\bar b=e^{2\sigma}b$ satisfy
\[
 \sigma(r_0)=c_0,
\]
\begin{equation}\label{eq-jet1}
 \sigma'(r_0)=-\frac{a'(r_0)}{2a(r_0)},
\end{equation}
\begin{equation}\label{eq-jet2}
 \sigma''(r_0)
 =\frac{2a(r_0)b(r_0)e^{2c_0}C_*-a''(r_0)
 +4a(r_0)\sigma'(r_0)^2}{2a(r_0)}.
\end{equation}
At $r_0$,
\begin{equation}\label{eq-norm-bara2}
 \bar a''=e^{2\sigma}\bigl(a''+2a\sigma''+4\sigma'a'+4a(\sigma')^2\bigr),
\end{equation}
and the criticality identity yields
\begin{equation}\label{eq-norm-bara2-placed}
 \bar a''=e^{2\sigma}\bigl(a''+2a\sigma''-4a(\sigma')^2\bigr).
\end{equation}
Consequently,
\begin{equation}\label{eq-local-twojet-optical-coefficient}
 \bar C_{\rm opt}(r_0)
 =\frac{\bar a''(r_0)}{2\bar a(r_0)\bar b(r_0)}
 =C_*.
\end{equation}
For a constant-speed angular geodesic with $c_\gamma=\breve\alpha(\gamma,\dot\gamma)>0$,
\begin{equation}\label{eq-twojet-affine-radial-mode}
 \delta\ddot r=\frac{\bar a''(r_0)}{2\bar b(r_0)}c_\gamma^2\delta r.
\end{equation}
For every nonzero angular pole $\breve v$ at $r_0$,
\begin{equation}\label{eq-local-twojet-radial-flag}
 \mathbf K_{\bar F}\bigl((r_0,u);(0,\breve v);(1,0)\bigr)=-C_*.
\end{equation}
When $C_*>0$, $\sqrt{C_*}$ is the real optical growth rate. When $C_*<0$, $\sqrt{-C_*}$ is the radial frequency. For $C_*=0$, the radial equation is degenerate and its solutions are affine.

Assume that an oriented closed angular geodesic is given. Then its lift at $r_0$ is an oriented closed constant-radius optical geodesic. If $a=r^2/A$ and $b=B/A$ for positive static coefficients $A,B$, its future null lift has periodic spatial projection and static time increment equal to its optical period.

Assume further that $A$ tends to one on a static spherical end, that $\sigma$ vanishes near that end, and that $S_*>0$ is a normalized shadow scale. Then
\begin{equation}\label{eq-jet0}
 c_0=\ln\!\left(\frac{S_*\sqrt{A(r_0)}}{r_0}\right)
 =\ln\!\left(\frac{S_*}{\sqrt{a(r_0)}}\right).
\end{equation}
\end{proposition}

\begin{proof}
We write
\[
 c_1:=-\frac{a'(r_0)}{2a(r_0)}
\]
and write $c_2$ for the expression in \eqref{eq-jet2}, with
$\sigma'(r_0)$ replaced by $c_1$. We choose a function $\chi\in C_c^\infty(J)$ that equals one on a neighborhood of $r_0$ and define
\[
 \sigma(r)=\chi(r)
 \left(c_0+c_1(r-r_0)+\frac{c_2}{2}(r-r_0)^2\right).
\]
As $\chi$ is constant near $r_0$,
\[
 \sigma(r_0)=c_0,\qquad
 \sigma'(r_0)=c_1,\qquad
 \sigma''(r_0)=c_2.
\]
The function has compact support. On every specified noncompact end with a
neighborhood disjoint from that support, $\sigma_\infty=0$.

Differentiation of $\bar a=e^{2\sigma}a$ gives
\[
 \bar a'=e^{2\sigma}(a'+2a\sigma')
\]
and \eqref{eq-norm-bara2}. At $r_0$,
\[
 a'(r_0)+2a(r_0)c_1=0,
\]
so $\bar a'(r_0)=0$. Substitution yields
\eqref{eq-norm-bara2-placed}. By the definition of $c_2$,
\[
 a''(r_0)+2a(r_0)c_2-4a(r_0)c_1^2
 =2a(r_0)b(r_0)e^{2c_0}C_*.
\]
Hence
\[
 \bar a''(r_0)=2a(r_0)b(r_0)e^{4c_0}C_*,
 \qquad
 \bar a(r_0)\bar b(r_0)=a(r_0)b(r_0)e^{4c_0}.
\]
Their quotient proves \eqref{eq-local-twojet-optical-coefficient}. Applying Theorem~\ref{thm-jacobi-ring} to $\bar a$ and $\bar b$, we obtain \eqref{eq-twojet-affine-radial-mode}. An oriented closed angular geodesic yields an oriented closed constant-radius optical geodesic. Under the static representative, its future null lift has periodic spatial projection and static time increment equal to the optical period. From Theorem~\ref{thm-flag}, we obtain
\[
 \mathbf K_{\bar F}
 =-\frac{\bar a''(r_0)}
 {2\bar a(r_0)\bar b(r_0)},
\]
which yields \eqref{eq-local-twojet-radial-flag}. The sign alternatives follow
from the optical-arclength equation
$\ddot\xi=\bar C_{\rm opt}\xi$.

In the spherical specialization, \eqref{eq-shadow-dressed} and
$\sigma_\infty=0$ give
\[
 S_*=e^{c_0}\frac{r_0}{\sqrt{A(r_0)}}
 =e^{c_0}\sqrt{a(r_0)}.
\]
Solving for $c_0$ proves \eqref{eq-jet0}.
\end{proof}

\begin{corollary}\label{cor-angular-finsler-placement}
Let
\[
 F^2=b(r)(v^1)^2+a(r)\breve\alpha(u,\breve v)^2,
 \qquad a(r)>0,\qquad b(r)>0,
\]
be defined on $\Cone$ over an open interval $I$. Assume that $(\breve M,\breve\alpha)$ admits an oriented closed geodesic $\gamma$. Fix $r_0\in I$, a compact interval $J\Subset I$ containing $r_0$ in its interior, and $c_0,C_*\in\RR$. Then there exists $\sigma\in C_c^\infty(J)$ such that $\bar F=e^{\sigma(r)}F$ has the constant-radius optical geodesic $(r_0,\gamma)$ and
\[
 \mathbf K_{\bar F}\bigl((r_0,u);(0,\dot\gamma);(1,0)\bigr)=-C_*,
 \qquad
 \ddot\xi=C_*\xi
\]
under optical-arclength normalization. If $a=r^2/A$ and $b=B/A$ for positive static coefficients $A,B$, the future null lift has periodic spatial projection and static time increment equal to the optical period. If $\breve\alpha$ is non-Riemannian, the angular Cartan tensor of $\bar F$ is nonzero wherever the angular Cartan tensor of $\breve\alpha$ is nonzero. If $\breve\alpha$ is nonreversible, the opposite orientation can have a different period while the radial coefficient remains $C_*$.
\end{corollary}

\begin{proof}
By Proposition~\ref{prop-local-conformal-twojet-flexibility}, we choose the compactly supported two-jet so that
$\bar a'(r_0)=0$ and $\bar C_{\rm opt}(r_0)=C_*$.  The angular geodesic therefore
lifts to a constant-radius optical geodesic, and
\eqref{eq-local-twojet-radial-flag} and \eqref{eq-local-twojet-optical-coefficient}
give the stated flag and Jacobi identities.  By Proposition~\ref{prop-static-representative}, we obtain the stated static light-ring interpretation. From Proposition~\ref{prop-non-riem-angular-sector} and the orientation dependence of Finsler length, we obtain the two assertions concerning the angular Cartan tensor and the opposite orientation of the given closed angular geodesic.
\end{proof}

\begin{remark}\label{rem-radial-sign-conventions}
The sign of $\bar a''(r_0)$ distinguishes exponential affine growth from oscillatory affine behavior.  A normalized optical-time exponent additionally contains the factor
$1/(\bar a(r_0)\bar b(r_0))$ as in \eqref{eq-lyapunov-optical-dressed}.
\end{remark}

\begin{example}\label{ex-concrete-placement-radial-behavior}
Take $b(r)\equiv1$ and $a(r)=1+r^{2}$ on $\RR$.  To place a critical radius at $r_0=1$, choose $\sigma'(1)=-1/2$.  With $\sigma(1)=0$, equation~\eqref{eq-norm-bara2-placed} gives
\[
 \bar a''(1)=4\sigma''(1).
\]
Thus $\sigma''(1)<0$ gives oscillatory radial behavior,
$\sigma''(1)=0$ gives the degenerate case, and
$\sigma''(1)>0$ gives exponential radial growth.
\end{example}

\begin{example}\label{ex-create-radial-mode}
Assume $a'(r)\neq0$ on $I$.  Fix $r_0\in I$ and impose
$\sigma'(r_0)=-a'(r_0)/(2a(r_0))$.  Then $r_0$ is a critical radius for the deformed structure.  Its affine radial coefficient changes sign when
\[
 \sigma''(r_0)
 =\frac{4a(r_0)(\sigma'(r_0))^2-a''(r_0)}{2a(r_0)}.
\]
Values on opposite sides of this threshold give oscillatory or exponential radial behavior.  A closed angular geodesic turns the constant-radius orbit into an oriented closed optical geodesic; the light-ring interpretation additionally requires the static representative convention.
\end{example}

\subsection{A Nonreversible Randers Example}
\label{sec-randers}

\noindent
We apply the local construction to a navigation Randers sphere. The two oriented equatorial geodesics have unequal periods, while the radial conformal data assign one common radial coefficient.

\begin{example}\label{ex-randers-sector}
The oriented equatorial orbits and unequal optical periods coincide with the analysis in \cite{Caponio2026}. We combine those angular periods with the locally assigned common radial coefficient. Let $\breve M=\mathbb S^2$ with round metric
\[
 h_{\mathbb S^2}=d\theta^2+\sin^2\theta\,d\phi^2,
\]
and choose the Killing wind
\[
 W=\omega\,\partial_\phi,
 \qquad 0<\omega<1.
\]
We write $\lambda(\theta)=1-\omega^2\sin^2\theta$. The Zermelo navigation norm is
\begin{equation}\label{eq-explicit-randers-norm}
 \breve\alpha_\omega(\theta,\phi;\dot\theta,\dot\phi)
 =\frac{
 \sqrt{\lambda(\dot\theta^2+\sin^2\theta\,\dot\phi^2)
 +\omega^2\sin^4\theta\,\dot\phi^2}
 -\omega\sin^2\theta\,\dot\phi}
 {\lambda}.
\end{equation}
It is a strongly convex nonreversible Randers metric of constant positive flag curvature, \cite{BaoRoblesShen2004}. On the equator $\theta=\pi/2$, formula
\eqref{eq-explicit-randers-norm} reduces to
\begin{equation*}
 \breve\alpha_\omega(0,\dot\phi)
 =
 \begin{cases}
 \dfrac{\dot\phi}{1+\omega},&\dot\phi>0,\\[2mm]
 \dfrac{-\dot\phi}{1-\omega},&\dot\phi<0.
 \end{cases}
\end{equation*}
Because navigation by a Killing wind sends round geodesics to Randers geodesics, the navigation images of the two unit-speed round equators are the oriented unit-speed
$\breve\alpha_\omega$-geodesics
\begin{align*}
 \gamma_+(s)&=\left(\frac\pi2,(1+\omega)s\right),
 &0\le s\le \frac{2\pi}{1+\omega},\\
 \gamma_-(s)&=\left(\frac\pi2,-(1-\omega)s\right),
 &0\le s\le \frac{2\pi}{1-\omega}.
\end{align*}
They trace the same equator with opposite orientations, and their Finsler lengths are
\begin{equation*}
 L_+=\frac{2\pi}{1+\omega},
 \qquad
 L_-=\frac{2\pi}{1-\omega}.
\end{equation*}
Thus reversal changes the oriented length.

On the conic domain of $(0,\infty)\times\mathbb S^2$, we consider
\[
 F^2=(v^1)^2+r^3\breve\alpha_\omega^2.
\]
Since $a'(r)=3r^2>0$, the original structure has no constant-radius angular critical radius. Choose an interval $J=(1/2,3/2)\Subset(0,\infty)$ and a compactly supported smooth function $\sigma\in C_c^\infty(J)$ that, in a neighborhood of $r=1$, has the given two-jet
\[
 \sigma(1)=0,
 \qquad
 \sigma'(1)=-\frac32,
 \qquad
 \sigma''(1)=-1.
\]
These three equalities assign only the local jet of $\sigma$ at $r=1$; away from this neighborhood the cutoff can be chosen arbitrarily inside $J$. Then $\bar a'(1)=0$ and
\[
 \bar a''(1)
 =6+2(-1)-4\left(\frac94\right)=-5.
\]
Also $\bar a(1)=\bar b(1)=1$. Consequently, both oriented equatorial geodesics lift at
$r=1$ to oriented closed constant-radius optical geodesics with oscillatory radial mode; after the static representative convention is imposed, their future null lifts have light-ring spatial projections. From Theorem~\ref{thm-flag-jacobi}, we obtain
\begin{equation*}
 \mathbf K_{\bar F}\bigl((1,u);(0,\dot\gamma_\pm);(1,0)\bigr)
 =\frac52,
\end{equation*}
and the optical-time radial equation has oscillation frequency $\sqrt{5/2}$ for both orientations. Their optical periods are direction dependent:
\begin{equation*}
 \mathcal T_+=\frac{2\pi}{1+\omega},
 \qquad
 \mathcal T_-=\frac{2\pi}{1-\omega},
 \qquad
 \frac{\mathcal T_-}{\mathcal T_+}=\frac{1+\omega}{1-\omega}>1.
\end{equation*}
The radial coefficient and radial-transverse flag curvature are common to both orientations, whereas the periods distinguish the nonreversible angular geometry. The angular velocities never vanish, so both lifted curves remain in $\Cone$. At $r=1$, the angular flag curvature is the angular Randers flag curvature because $\bar a'(1)=0$ and $\bar a(1)=1$.
\end{example}

\subsection{Curvature Identities and Riemannian Specializations}\label{subsec-curvature-checks}

We verify the spray-kernel identity, metric self-adjointness, Ricci trace, and Riemannian reduction of the curvature block.

When the angular norm is Riemannian, we write $a=f^2$ and $b=1$. Then
\[
 P=ff'',\qquad Q=(f')^2,
\]
so the radial and angular sectional curvatures reduce to
\[
 -\frac{P}{a}=-\frac{f''}{f},
 \qquad
 \frac{\breve\kappa-Q}{a}=\frac{\breve\kappa-(f')^2}{f^2}.
\]
Finally, tracing the angular block yields
\[
 R^a{}_a
 =\breve\Ric-(n-1)\frac baP(v^1)^2-(n-2)Q\breve\alpha^2,
\]
and addition of $R^1{}_1=-P\breve\alpha^2$ reproduces
\eqref{eq-Ricci-decomp}.

When the angular norm is Riemannian, write $h=f^2$. Then,  we have
\[
 \mathbf K_{\rm rad}=-\frac{f''}{f}
 =-\frac{h''}{2h}+\frac{(h')^2}{4h^2},
 \qquad
 \mathbf K_{\rm ang}=\frac{\breve\kappa-(f')^2}{f^2}
 =\frac{\breve\kappa}{h}-\frac{(h')^2}{4h^2}.
\]
These expressions coincide with the Riemannian specialization of the flag-curvature formulas established in this section. For $n\ge3$, the Einstein reduction studied in the Finsler warped-product literature, \cite{ChenShenZhao2018} gives
\[
 \breve\Ric=(n-2)(Q-P)\breve\alpha^2,
 \qquad
 \lambda=-\frac{P}{a},
\]
when $\Ric_F=(n-1)\lambda F^2$; the common value $Q-P$ is constant on each connected component. These cited relations verify the explicit conic curvature block.

\subsection{Spherical Riemannian Shadow Scale}\label{subsec-shadow-scale}

\noindent
We consider a static, spherically symmetric exterior with a normalized time coordinate. The angular optical coefficient determines the critical impact parameter, and a radial conformal factor changes this quantity through the difference between its value at the critical radius and its limiting value.

Consider a static spherical exterior
\[
 ds^{2}=-A(r)dt^{2}+B(r)dr^{2}+r^{2}d\Omega^{2},
 \qquad A(r)>0,\quad B(r)>0,
\]
with $A\to1$ in the asymptotic region. For a spherical critical radius with exponential radial sign, the established critical impact scale is \cite{ClaudelVirbhadraEllis2001,Qiao2022}
\begin{equation*}
 b_{\rm sh}(r_0)=\sqrt{a(r_0)}=\frac{r_0}{\sqrt{A(r_0)}}.
\end{equation*}
For the optical deformation $\bar F=e^{\sigma(r)}F$, assume that $\sigma_\infty=\lim_{r\to\infty}\sigma(r)$ exists and that the deformed optical metric is connected with a static representative whose time coordinate is normalized in the asymptotic region. Then
\begin{equation}\label{eq-shadow-dressed}
 \bar b_{\rm sh}^{\rm norm}(r_0)
 =e^{\sigma(r_0)-\sigma_\infty}\frac{r_0}{\sqrt{A(r_0)}}.
\end{equation}
Compact support gives $\sigma_\infty=0$. This scalar relation is confined to the Riemannian spherical specialization. A nonreversible angular Finsler sector has direction-dependent critical impact data, \cite{Caponio2026}.

\begin{remark}
The operation $F\mapsto e^\sigma F$ is a deformation of the Fermat metric.  It changes its geodesics unless $\sigma$ is constant.  This operation is distinct from multiplication of the full spacetime metric by a conformal factor, which preserves unparameterized null geodesics.
\end{remark}

\begin{remark}
For the Schwarzschild value $b_{\rm sh}=3\sqrt3\,M$, formula~\eqref{eq-shadow-dressed} gives
\[
 \bar b_{\rm sh}^{\rm norm}
 =e^{\sigma(r_0)-\sigma_\infty}3\sqrt3\,M.
\]
This scaling statement is meaningful only after the asymptotic time normalization has been fixed.
\end{remark}

\begin{example}\label{ex-shadow-design}
For a chosen exterior radius $r_0$ and target shadow scale $S_*$, equation~\eqref{eq-jet0} determines the conformal value. Equation~\eqref{eq-jet1} places the critical radius. After these assignments, equation~\eqref{eq-jet2} determines $\sigma''(r_0)$ for a selected signed radial coefficient. These assignments are successive.
\end{example}

\begin{example}\label{ex-quadratic-jet-design}
Let $S_*>0$ and a signed radial coefficient $C_*\in\RR$ be given at $r_0$. We write
\begin{align*}
 c_0&=\ln\!\left(\frac{S_*}{\sqrt{a(r_0)}}\right),\\
 c_1&=-\frac{a'(r_0)}{2a(r_0)},\\
 c_2&=\frac{2a(r_0)b(r_0)e^{2c_0}C_*
 -a''(r_0)+4a(r_0)c_1^2}{2a(r_0)}.
\end{align*}
We choose a cut-off function $\chi$ that equals one near $r_0$ and has compact support in the exterior, and define
\[
 \sigma(r)=\chi(r)
 \left(c_0+c_1(r-r_0)+\frac{c_2}{2}(r-r_0)^2\right).
\]
The resulting deformation has $\sigma_\infty=0$ and realizes the three specified local quantities.
\end{example}
\begin{remark}\label{rem-design-recipe}
The conformal data are assigned in order: $\sigma(r_0)-\sigma_\infty$ determines the normalized shadow scale in the spherical exterior, $\sigma'(r_0)$ determines the critical radius, and $\sigma''(r_0)$ determines the normalized radial coefficient once the value and slope have been chosen.
\end{remark}
\subsection{Spherical Optical Examples}\label{subsec-riemannian-checks}

We compare the conic formulas with the standard Riemannian warped-product identities and evaluate the radial coefficients in Schwarzschild, Reissner\textendash{}Nordstr\"om, and Kottler geometries.

Let $(\breve M,\breve\alpha)$ be Riemannian with constant sectional curvature $\breve\kappa$. 
\[
 \mathbf K_{\rm ang}=\frac{\breve\kappa}{h}-\frac{(h')^2}{4h^2},
 \qquad
 \mathbf K_{\rm rad}=-\frac{h''}{2h}+\frac{(h')^2}{4h^2},
\]
and
\[
 \Ric_F=-(n-1)\left(\frac{h''}{2h}-\frac{(h')^2}{4h^2}\right)(v^1)^2
 +\left((n-2)\breve\kappa-\frac{h''}{2}-\frac{(n-3)(h')^2}{4h}\right)\breve\alpha^2.
\]
Substitution of a Riemannian angular norm into the conic curvature and flag formulas reproduces these identities. We evaluate the radial coefficients in three standard static geometries.

\begin{example}\label{ex-schwarzschild}
For the Schwarzschild exterior $r>2M$,
\[
 A(r)=1-\frac{2M}{r},\quad B(r)=A(r)^{-1},\qquad
 a(r)=\frac{r^3}{r-2M},\quad b(r)=\frac{r^2}{(r-2M)^2}.
\]
Direct differentiation yields
\[
 a'(r)=\frac{2r^2(r-3M)}{(r-2M)^2},
 \qquad
 a''(r)=\frac{2r(r^2-6Mr+12M^2)}{(r-2M)^3}.
\]
The exterior critical radius is $r_0=3M$, with
\[
 a(3M)=27M^2,\qquad b(3M)=9,\qquad a''(3M)=18.
\]
The affine growth coefficient and normalized optical-time exponent are
\[
 \kappa_{\rm aff}(3M)=\sqrt{\ell},
 \qquad
 \gamma_{\rm opt}(3M)=\frac{1}{3\sqrt3\,M}.
\]
The asymptotically normalized shadow radius is
\[
 b_{\rm sh}(3M)=3\sqrt3\,M.
\]
For the round angular factor, the angular and radial sectional curvatures at the critical radius are
\[
 \mathbf K_{\rm ang}=\frac{1}{27M^2},
 \qquad
 \mathbf K_{\rm rad}=-\frac{1}{27M^2}.
\]
Thus the radial optical curvature has the sign associated with exponential radial growth.
\end{example}

\begin{example}\label{ex-RN}
For the Reissner\textendash Nordstr\"om geometry,
\[
 A(r)=1-\frac{2M}{r}+\frac{Q^2}{r^2},
 \qquad B(r)=A(r)^{-1},
\]
one has
\[
 a(r)=\frac{r^4}{r^2-2Mr+Q^2},
 \qquad
 b(r)=\frac{r^4}{(r^2-2Mr+Q^2)^2},
\]
and
\[
 a'(r)
 =\frac{2r^3(r^2-3Mr+2Q^2)}
 {(r^2-2Mr+Q^2)^2}.
\]
The algebraic critical radii are
\begin{equation*}
 r_{\rm ph}^{\pm}
 =\frac{3M\pm\sqrt{9M^2-8Q^2}}{2},
\end{equation*}
whenever $9M^2\ge8Q^2$.  At a nonzero root $r_0$,
\[
 \frac{a''(r_0)}{2b(r_0)}
 =\frac{2r_0-3M}{r_0}.
\]
Hence $r_{\rm ph}^{+}$ has the exponential radial sign and
\[
 \kappa_{\rm aff}(r_{\rm ph}^{+})
 =\sqrt{\frac{2r_{\rm ph}^{+}-3M}{r_{\rm ph}^{+}}\,\ell},
\]
while
\[
 \gamma_{\rm opt}(r_{\rm ph}^{+})
 =\sqrt{
 \frac{(2r_{\rm ph}^{+}-3M)
 (r_{\rm ph}^{+\,2}-2Mr_{\rm ph}^{+}+Q^2)}
 {r_{\rm ph}^{+\,5}}}.
\]
For $0<|Q|<M$, the event horizons are
\[
 r_{\rm h}^{\pm}=M\pm\sqrt{M^2-Q^2}.
\]
The radius $r_{\rm ph}^{+}$ lies in the static exterior $r>r_{\rm h}^{+}$. The algebraic radius $r_{\rm ph}^{-}$ lies between the two horizons; its quadratic-variation sign refers to the analytic continuation of the coefficient. The extremal and horizonless cases require separate domain conditions.
\end{example}

\begin{example}\label{ex-kottler}
For the Kottler static region,
\[
 A(r)=1-\frac{2M}{r}-\frac{\Lambda r^2}{3},
 \qquad B(r)=A(r)^{-1},
 \qquad A(r)>0,
\]
the criticality equation $2A=rA'$ gives $r_0=3M$.  At this radius,
\[
 A(3M)=\frac13-3\Lambda M^2,
 \qquad
 \frac{a''(3M)}{2b(3M)}=1.
\]
Thus
\[
 \kappa_{\rm aff}(3M)=\sqrt{\ell},
 \qquad
 \sqrt{a(3M)}
 =\frac{3M}{\sqrt{1/3-3\Lambda M^2}}.
\]
This expression is a local optical scale. Since Kottler geometry has no asymptotically flat end, an observer-defined shadow requires a specified observer position and normalization.
\end{example}

\subsection{Gauss\textendash{}Bonnet Calculation on the Rotational Optical Surface}\label{subsec-gauss-bonnet}

\noindent
We apply the Gauss\textendash{}Bonnet theorem to a rotational optical annulus and compute the curvature integral together with both oriented boundary contributions.

For the two-dimensional rotational optical metric
\(
 dl_{\rm opt}^{2}=b(r)dr^{2}+a(r)d\phi^{2},
 \qquad a(r)>0,\quad b(r)>0,
\)
we have
\begin{equation}\label{eq-gaussian}
 K_{\rm opt}
 =-\frac{1}{\sqrt{ab}}\frac{d}{dr}\left(\frac{a'}{2\sqrt{ab}}\right)
 =-\frac{a''}{2ab}+\frac{a'b'}{4ab^{2}}+\frac{(a')^{2}}{4a^{2}b}.
\end{equation}
At a critical circle, condition~\eqref{eq-angular-critical-known} yields
\begin{equation}\label{eq-gaussian-at-ring}
 K_{\rm opt}(r_0)=-\frac{a''(r_0)}{2a(r_0)b(r_0)},
 \qquad
 C_{\rm opt}(r_0):=-K_{\rm opt}(r_0).
\end{equation}
For $C_{\rm opt}>0$, we have $\gamma_{\rm opt}=\sqrt{C_{\rm opt}}$. For $C_{\rm opt}<0$, solutions of the radial equation are oscillatory.

For the equatorial annulus
\[
 \Omega=[r_-,r_+]\times S^1,
\]
the Gauss-Bonnet theorem yields $\chi(\Omega)=0$. Direct integration of \eqref{eq-gaussian} yields
\begin{equation}\label{eq-GB-optical}
 \int_{\Omega}K_{\rm opt}\,dA_{\rm opt}
 =-\pi\left[\frac{a'}{\sqrt{ab}}\right]_{r_-}^{r_+},
\end{equation}
and the two oriented boundary circles give
\begin{equation}\label{eq-lyapunov-Kopt}
 \int_{\partial\Omega}\kappa_g\,ds
 =\pi\left[\frac{a'}{\sqrt{ab}}\right]_{r_-}^{r_+}.
\end{equation}
Consequently,
\[
 \int_{\Omega}K_{\rm opt}\,dA_{\rm opt}
 +\int_{\partial\Omega}\kappa_g\,ds=0.
\]
This identity is consistent with the local radial coefficient.

\begin{remark}\label{rem-GB-physics}
Equation~\eqref{eq-gaussian-at-ring} identifies the Gaussian curvature with the negative radial coefficient at a critical circle. Equations~\eqref{eq-GB-optical} and \eqref{eq-lyapunov-Kopt} show that the two oriented boundary contributions cancel the annular curvature integral.
\end{remark}

\subsection*{Declarations}

\paragraph{Funding.}
Not applicable.

\paragraph{Conflict of interest.}
The authors declare no conflict of interest.

\paragraph{Ethics approval.}
This study does not involve human participants or animals.

\paragraph{Data availability.}
No dataset is generated or analyzed in this work.

\paragraph{Author contributions.}
The authors contributed to the conception, formulation, analysis, writing, and revision of the text.

\end{document}